\documentclass[12pt]{amsart}

\usepackage{amsmath}
\usepackage{amssymb}
\usepackage{fancybox}
\usepackage{eucal}
\usepackage{amsthm}
\usepackage{amscd}
\usepackage{wasysym}
\usepackage{url}
\usepackage{MnSymbol} 

\usepackage[all,cmtip]{xy}

\usepackage{enumitem}

\usepackage{hyperref}

\newtheorem{thm}{Theorem}[section]
\newtheorem{prop}[thm]{Proposition}

\newtheorem{lemma}[thm]{Lemma}

\newtheorem{cor}[thm]{Corollary}

\newtheorem{notation}[thm]{Notation}

\newtheorem{definitiontemp}[thm]{Definition}
\newenvironment{defn}{\begin{definitiontemp} 
\normalfont}{\end{definitiontemp}}

\newenvironment{pf}{\begin{trivlist}\item[\hskip\labelsep
{\it Proof.}]}{\end{trivlist}}

\newcommand{\comment}[1]{}

\newcommand{\ACF}[1]{\ensuremath{\textbf{ACF}_{#1}}}
\newcommand{\DCF}[1]{\ensuremath{\textbf{DCF}_{#1}}}

\newcommand{\set}[2]{\ensuremath{ \{ #1 : #2 \} }}

\newcommand{\ord}[1]{\ensuremath{ \text{ord} (#1 )}}
\newcommand{\order}[2]{\ensuremath{ \text{ord}_{#2} (#1 )}}
\newcommand{\spec}[1]{\ensuremath{ \text{Spec} (#1 )}}

\renewcommand{\deg}[1]{\ensuremath{\text{deg}(#1)}}

\newcommand{\Z}{\mathbb{Z}}

\newcommand{\Q}{\mathbb{Q}}
\newcommand{\R}{\mathcal{R}}
\newcommand{\C}{\mathcal{C}}
\newcommand{\A}{\mathcal{A}}
\newcommand{\B}{\mathcal{B}}
\newcommand{\D}{\mathcal{D}}

\newcommand{\Gtilde}{\widetilde{G}}
\newcommand{\Htilde}{\widetilde{H}}

\newcommand{\ytilde}{\tilde{y}}

\newcommand{\M}{\mathfrak{M}}

\renewcommand{\S}{\mathcal{S}}

\newcommand{\Xvec}{\vec{X}}
\newcommand{\yvec}{\vec{y}}

\newcommand{\del}{\delta}

\newcommand{\Ky}{K\la y\ra}
\newcommand{\KY}{K\{ Y\}}
\newcommand{\Trans}{\text{Trans}}

\newcommand{\proves}{\ensuremath{\vdash}}

\newcommand{\la}{\langle}
\newcommand{\ra}{\rangle}

\def\converges{\!\downarrow}
\newcommand{\at}{\char'100}

\newcommand{\dom}[1]{\text{dom}(#1)}
\newcommand{\rg}[1]{\text{range}(#1)}

\newcommand{\bfb}{\boldsymbol{b}}
\newcommand{\bfc}{\boldsymbol{c}}
\newcommand{\bfd}{\boldsymbol{d}}
\newcommand{\Qhat}{\widehat{\mathbb{Q}}}

\newcommand{\Khat}{\widehat{K}}
\newcommand{\Lhat}{\widehat{L}}
\def\bfz{\boldsymbol{0}}
\def\s01{\ensuremath{\Sigma^0_1}}
\def\d02{\ensuremath{\Delta^0_2}}
\def\phi{\varphi}

\def\res{\!\!\upharpoonright\!}

\begin{document}

\title[Turing Degree Spectra of Differentially Closed Fields]{Turing Degree Spectra of\\Differentially Closed Fields}
\author[D. Marker]{David Marker}
\address{Dept.\ of Mathematics, Statistics, \& Computer Science, 
University of Illinois at Chicago,
Chicago, IL, U.S.A.}
\email{marker@math.uic.edu}

\author[R. Miller]{Russell Miller}
\address{Dept.\ of Mathematics, Queens College, \& Ph.D.\ Programs in Mathematics \& Computer Science,
Graduate Center, City University of New York,  USA}
\email{Russell.Miller@qc.cuny.edu}
\urladdr{\href{http://qcpages.qc.cuny.edu/~rmiller/}{http://qcpages.qc.cuny.edu/$\sim$rmiller/}}

\thanks{
The second author was
partially supported by Grants \# DMS -- 1001306 \& 1362206 from
the National Science Foundation, and by several grants from
the PSC-CUNY Research Award Program.  This work was initiated
at a workshop held at the American Institute of Mathematics
in August 2013, where the question of noncomputable
differentially closed fields was raised by Wesley Calvert.
The authors appreciate the support of A.I.M., and also
thank Uri Andrews and Hans Schoutens for useful conversations.\\
\indent
The authors and all who read the proof of Theorem \ref{thm:low} owe a debt of gratitude
to the anonymous but valiant referee for a number of
suggestions which substantially simplified that proof.}

\begin{abstract}
The degree spectrum of a countable structure is the set of all
Turing degrees of presentations of that structure.  We show that
every nonlow Turing degree lies in the spectrum of some 
differentially closed field (of characteristic $0$, with a single derivation)
whose spectrum does not contain the computable degree $\bfz$.
Indeed, this is an equivalence, for we also show that 
if this spectrum contained a low degree,
then it would contain the degree $\bfz$.
From these results we conclude
that the spectra of differentially closed fields of
characteristic $0$ are exactly the jump-preimages of spectra
of automorphically nontrivial graphs.
\end{abstract}

\maketitle

\section{Introduction}
\label{sec:intro}

Differential fields arose originally in work of Ritt examining algebraic
differential equations on manifolds over the complex numbers.
Subsequent work by Ritt, Kolchin and others brought this study
into the realm of algebra, where numerous parallels appeared
with algebraic geometry.  The topic first intersected with model
theory in the mid-twentieth century, in work of
Abraham Robinson, and logicians soon discovered
the theories of differential fields and of differentially closed fields
to have properties which had been considered in the abstract,
but had not previously been known to hold for any everyday theories
in mathematics.  It was the model theorists who provided
the definitive resolution to the question of differential closure,
several variations of which had previously been
developed in differential algebra.
In 1974, Harrington proved the existence of computable
differentially closed fields, making the notion more concrete,
although our grasp of this topic remains more tenuous
than our understanding of algebraic closures in field theory.

In this article, we offer an analysis of the complexity
of countable models of \DCF0,
the theory of differentially closed fields of characteristic $0$.
This work requires a solid background in differential algebra,
in model theory, and in effective mathematics.  Ultimately
we will characterize the spectra of countable models of \DCF0
as exactly the preimages, under the jump operation, of spectra
of automorphically nontrivial countable graphs; or, equivalently,
as exactly those spectra of such graphs which are closed under
a simple equivalence relation on Turing degrees.  To do so, we show
that spectra of differentially closed fields have certain complexity
properties, which are not known to hold of any other standard
class of mathematical structures:  every low differentially closed field of
characteristic $0$ is isomorphic to a computable one, whereas every
nonlow degree computes a differentially closed field which has no
computable copy.  Indeed we will present a substantial class of
fairly complex spectra that can all be realized by models of \DCF0,
including spectra with arbitrary proper $\alpha$-th jump degrees,
for every computable nonzero ordinal $\alpha$.  To explain what
these results mean, we begin immediately with the necessary background.
For supplemental information on computability theory,
\cite{S87} is a standard source, while for more detail about
model theory and differential fields, we suggest \cite{M06},
\cite{P02}, or the earlier \cite{S72}.

\subsection{Background in Differential Algebra}
\label{subsec:diffalg}

A differential ring is a ring with a \emph{differential
operator}, or \emph{derivation}, on its elements.  If the ring
is a field, we call it a \emph{differential field}.  The differential
operator $\del$ is required to preserve addition and to satisfy the
familiar \emph{Leibniz Rule}:  $\del (x\cdot y)=(x\cdot\del y)+(y\cdot\del x)$.
Examples include the field $\Q(x)$ of rational functions over $\Q$
in a single variable $x$, with the usual differentiation $\frac{d}{dx}$,
or the field $\Q(t,\del t,\del^2 t,\ldots)$, with $\del$ acting 
as suggested by the notation.  In these examples, $\Q$ may be replaced
by another differential field $K$, with the derivation $\del$ on $K$
likewise extended to all of $K(x)$ or $K(t,\del t,\ldots)$.
(The only possible derivation on $\Q$ maps all rationals to $0$.
In general, the \emph{constants} of a differential field $K$ are those
$x\in K$ with $\del x=0$, and they form the \emph{constant
subfield} $C_K$ of $K$.)
We use angle brackets and write $K\la y_i~:~i\in I\ra$ for the smallest
differential subfield (of a given extension of $K$) containing
all the elements $y_i$; this is well-defined, and this subfield is said
to be \emph{generated} as a differential field by $\{ y_i~|~i\in I\}$.
Of course, the field generated by these same elements $\{ y_i~|~i\in I\}$ may well
be a proper subfield of this:  in the examples above, $\Q\la x\ra=\Q(x)$,
but $\Q(t)\subsetneq\Q\la t\ra=\Q(t,\del t,\ldots)$.
Differentiation of rational functions turns out to follow the usual
quotient rule, noting that $\del$ may map coefficients
in a nonconstant ground field $K$ to elements other than $0$.

For the purposes of this article, we restrict ourselves to characteristic $0$
and to \emph{ordinary} differential rings and fields, i.e., those with
only one derivation.  \emph{Partial} differential rings,
with more differential operators, exist and have natural
examples, as do differential rings of positive characteristic,
but considering either would expand this article well beyond
the scope we intend.

For a differential ring $K$ with derivation $\del$, $K\{ Y\}$ denotes
the ring of all \emph{differential polynomials} over $K$; it may be
viewed as the ring of algebraic polynomials $K[Y,\del Y,\del^2Y,\ldots]$,
with $Y$ and all its derivatives treated as separate variables.  We then define
$K\{ Y_0,\ldots,Y_{n+1}\} = (K\{ Y_0,\ldots,Y_n\})\{ Y_{n+1}\}$.
One sometimes differentiates a differential polynomial, treating
each $\del^{n+1}Y_i$ as the derivative of $\del^nY_i$.  
With only one derivation in the language, we often write $Y'$
for $\del Y$, or $Y^{(r)}$ for $\del^rY$.
\comment{
If $K$ has several derivations, we allow them all to appear
in the monomials, but usually they are assumed to commute
with each other in $K$, and if so, then we take them to commute
with each other in the monomials as well, so that $(\del^3_1\del_2Y)^2Y^5$
is the same monomial as $(\del_2\del_1^3Y)^2Y^5$, for instance.
}

The \emph{order} of a nonzero differential polynomial
$q\in K\{ Y\}$ is the greatest $r$ such that the $r$-th derivative
$Y^{(r)}$ appears nontrivially in $q$.  Equivalently, it is the least $r$
such that $q\in K[Y,Y',\ldots,Y^{(r)}]$.  Having order $0$ means that $q$
is an algebraic polynomial in $Y$ of degree $>0$; nonzero elements
of $K$ within $\KY$ are said to have order $-1$, and in this article,
the order of the zero polynomial is taken to be $+\infty$.  Each nonzero polynomial 
in $K\{ Y\}$ also has a \emph{rank} in $Y$.  For two such polynomials,
the one with lesser order has lesser rank.  If they have the same order $r$, then
the one of lower degree in $Y^{(r)}$ has lesser rank.  
Having the same order $r$ and the same degree in $Y^{(r)}$
is sufficient to allow us to reduce one of them, modulo the other,
to a polynomial of lower degree in $Y^{(r)}$, and hence of lower rank:
just take an appropriate $K$-linear combination of the two.
So, for our purposes, the rank in $Y$ is simply given by the order $r$
and the degree of $Y^{(r)}$.  Therefore, our ranks of nonzero differential
polynomials will be ordinals in $\omega^2$.

Our convention in this article is that the zero polynomial
has order $+\infty$.  Thus, for every element $x$ in any differential field
extension of $K$, the \emph{minimal differential polynomial} of $x$ over $K$
is defined (up to a nonzero scalar from $K$) as the differential polynomial
$q$ in $K\{ Y\}$ of least rank for which $x$ is a zero (i.e., $q(x)=0$).
In particular, the zero polynomial is considered to be the minimal
differential polynomial of an element differentially transcendental over $K$
(such as $t$ in $\Q\la t\ra$ above);
this is simply for notational convenience.

The \emph{differential closure} $\Khat$ of a differential field $K$
is the prime model of the theory $\DCF0\cup\Delta(K)$,
the union of the atomic diagram $\Delta(K)$ of $K$
with the (complete) theory \DCF0.
This theory was effectively axiomatized by Blum:
her axioms for a differentially closed field $F$
include the axioms for differential fields of characteristic $0$
and state, for each pair $(p,q)$ of differential polynomials
with arbitrary coefficients from $F$ and
with $\ord{p}>\ord{q}$, that $F$ must contain
an element $x$ with $p(x)=0\neq q(x)$.
(By our convention on ranks, $\ord{p}>\ord{q}$ ensures
that $q$ is not the zero polynomial.  However,
$q$ may equal $1$, and so $F$ must be algebraically closed.
Notice here that, for all fields, model-theoretic algebraic
closure implies field-theoretic algebraic closure.)
Blum's proofs appear in \cite{B68,B77}, and a summary
of all these results can be found in \cite{M06}.

Abraham Robinson
showed that \DCF0 has quantifier elimination.  Blum's computable
axiomatization makes \DCF0 decidable, hence makes the quantifier
elimination effective, both of which are particularly important for work
involving computable-model-theoretic questions about \DCF0.
Every definable set in a computable model of \DCF0 must now be
decidable, and, given the original defining formula of the set,
we can effectively find an equivalent quantifier-free formula, thereby
passing uniformly to the decision procedure for the set.
(Of course, this applies only to finitary defining formulas,
not to computable infinitary formulas.)

Blum proved \DCF0 to be $\omega$-stable, and existing
results of Morley then established that the theory $\DCF0\cup\Delta(K)$
always has a prime model, i.e., every differential field $K$ has
a differential closure.  Subsequently, Shelah proved
that, as the prime model extension of an $\omega$-stable theory, 
the differential closure $\Khat$ of $K$ is unique and realizes exactly those types
principal over $K$.  Each principal $1$-type has as generator a formula of the form
$p(Y)=0\neq q(Y)$, where $(p,q)\in (K\{ Y\})^2$ is a \emph{constrained pair}.
By definition, this means that $p(Y)$ is a monic, algebraically irreducible
polynomial in $K\{ Y\}$, that $q$ has strictly lower rank in $Y$ than $p$
does, and that, in $\Khat$ (and hence in every differential field extension of $K$),
every $y$ satisfying $p(y)=0\neq q(y)$ has minimal differential polynomial $p$
over $K$.  (A fuller definition appears in \cite[Defn.\ 4.3]{MOT14}.)
Hence the elements satisfying the generating formula form an orbit
under the action of those automorphisms of $\Khat$ that fix $K$ pointwise.
For a pair $(p,q)$ to be constrained is a $\Pi^K_1$ property, and there exist
computable differential fields $K$ for which it is $\Pi_1$-complete.
(This can happen even for a constant field $K$, such as the field
$\Q[\sqrt{p_n}~:~n\in\emptyset']$; see \cite{M08}.)
If such a $q$ exists, then $p$ is said to be \emph{constrainable};
clearly this property is $\Sigma^K_2$.
Not all monic irreducible polynomials in $K\{ Y\}$ are constrainable:
for example, $\del Y$ is not.  More generally, no $p$
in the image of $K\{ Y\}$ under $\del$ is constrainable, and
certain polynomials $p$ outside this image are also known to be unconstrainable.
In fact, constrainability has been shown in \cite{MOT14}
to be $\Sigma^0_2$-complete for certain computable differential fields $K$.
The exact complexity of constrainability over the constant differential
field $\Q$ is unknown:  it might even be decidable.
We note that $p$ is constrainable over $K$ if
\label{mdpconstrainable}
and only if some $y\in\Khat$ has minimal differential polynomial $p$ over $K$.
(This equivalence will be extremely useful in the $\S_m$-substages of the
construction for Theorem \ref{thm:low}.)  The equivalent condition proves
again that constrainability is $\Sigma^K_2$, provided that there exists a
$K$-computable presentation of $\Khat$, which we get from a theorem of Harrington.

\begin{thm}[Harrington; \cite{H74}, Corollary 3]
\label{thm:Harrington}
For every computable differential field $K$, there exists a
computable differential field $L$ and a computable differential
field homomorphism $g:K\to L$ such that $L$ is a differential closure of
the image $g(K)$.  Moreover, indices for $g$ and $L$ may be found
uniformly in an index for $K$.
\end{thm}
So this $L$ is in fact a differential closure of $K$ -- or at least,
of the image $g(K)$, which is computably isomorphic to $K$ via $g$.
In \cite{R60}, Rabin proved the original analogue of this theorem for 
fields and their algebraic closures.
We note that the exposition in \cite{H74} does not consider
uniformity of the procedure it describes, but a close reading of the proof
there indicates that the algorithm giving $g$ and $L$ is indeed
uniform in an index for $K$.  In particular, the following lemma
is proven simply by uniformizing the proof of \cite[2(b), Lemma 2]{H74}
and noting that the argument in the ensuing section 2(c) is uniform.
\begin{lemma}
\label{lemma:typefct}
There exists a single computable function $\gamma$, the \emph{type function}
for \DCF0, such that, for every computable differential field $F$
of characteristic $0$, every index $e$ for the atomic diagram $\Delta(F)$ of $F$,
and every irreducible differential polynomial $f\in F\{ X\}$,
$\gamma(e,f)$ is an index of the characteristic function $\phi_{\gamma(e,f)}$
of a $1$-type $\Gamma(x)$ that is complete and principal over
$\DCF0\cup\Delta(F)$ and contains the formula $f(X)=0$.
\qed\end{lemma}
This type function will enable us to extend individual formulas
$f(X)=0$ uniformly to principal $1$-types over
differential fields we have already built.  However, while the type
$\Gamma(x)$ given by the type function will always be principal,
the lemma does not promise to identify any specific formula
as a generator of the type.  The characteristic function merely
decides which formulas belong to the type and which do not:
at some point it will come across a generating formula and include it,
but having done so, it will simply continue including and excluding
other formulas, although from then on the type is in fact completely determined.

\subsection{Background in Model Theory}
\label{subsec:modeltheory}

Proposition \ref{prop:DCFtoGraph} will require some background beyond
Subsection \ref{subsec:diffalg}, which we provide here, referring
the reader to \cite{M06} and \cite{pillay} 
for details and further references regarding these results.
Model theorists have made dramatic inroads in the study of differential
fields and \DCF0; here we restrict ourselves to describing the results
necessary to prove Proposition \ref{prop:DCFtoGraph}, without giving complete
definitions of all the relevant concepts.

Let $K$ be a differentially closed field, with subfield $C_K$ of constants.
For $a\in K\setminus C_K$, consider the elliptic curve $E_a$ given by 
$$y^2=x(x-1)(x-a).$$
Let $E^\sharp_a$ be the Kolchin closure of the set of all torsion points
in the usual group structure on $E_a$.  (The Kolchin topology
is the differential analogue of the Zariski topology.)  The set $E^\sharp_a$
is known as the \emph{Manin kernel} of this abelian variety, as it is the
kernel of a certain homomorphism of differential algebraic groups.
One construction of Manin kernels appears in \cite{marker-manin}.
In the proof of Proposition \ref{prop:DCFtoGraph} we will use Manin
kernels $E^\sharp_{a_ma_n}$, meaning $E^\sharp_{a}$ as above
with $a=a_m+a_n$.

\begin{thm}
\label{thm:Maninkernel}
The family $\set{E^\sharp_a}{a^\prime\ne 0}$ is definable.
Indeed, it can be defined uniformly in each $a$ with
$a^\prime=a^3-a^2\neq 0$, by a quantifier-free formula.
\end{thm}

The definability is claimed in \cite{hi} but done more clearly in \cite[Sec.\ 2.4]{np}.
Of course, quantifier elimination for \DCF0 allows us to take the definition to be quantifier-free.
The condition $a^\prime=a^3-a^2$ will be relevant below.

\begin{thm}
\label{thm:Maninprops}
If $a^\prime\ne 0$, then $E^\sharp_a$ is strongly minimal and locally modular.
Moreover,
$E^\sharp_a$ and $E^\sharp_b$ are non-orthogonal if and only if $E_a$ and $E_b$ are isogenous.
In particular if $a$ and $b$ are algebraically independent over $\Q$, then $E^\sharp_a$
and $E^\sharp_b$ are orthogonal.
\end{thm}

These results are due to Hrushovski and Sokolovi\'c \cite{hs}, whose manuscript was never
published.   A proof  of the first fact is given in Section 5 of \cite{marker-manin}, and proofs of
both results appear in Section 4 of \cite{pillay}.

\begin{cor}
For every element $(b_0,b_1)$ of $E^\sharp_a$ in the differential closure of $\Q\la a\ra $,
both $b_0$ and $b_1$ are algebraic over $\Q\la a\ra $.
\end{cor}

\begin{pf}
Let $\psi(b_0,b_1)$ be the formula over $\Q\la a\ra $ isolating the type of $(b_0,b_1)$.
If $\psi$ defined an infinite subset of $E^\sharp_a$, then it would contain a torsion point.
But if $\psi$ contains an $n$-torsion point, every point in $\psi$ would be an $n$-torsion point,
yet there are only $n^2$ $n$-torsion points in $E_a$, a contradiction.
Thus $\psi(b_0,b_1)$ defines a finite set, so this pair is model-theoretically algebraic over $a$,
hence lies in the field-theoretic algebraic closure of $\Q\la a\ra $.
\qed\end{pf}

\begin{lemma}
\label{lemma:sms}
Let $X$ and $Y$ be strongly minimal sets defined over a differentially closed field $K$.
If  $X$ and $Y$ are orthogonal, then for any new element $x\in X$
the differential closure of $K\la x\ra $ contains no new elements of $Y$.
\end{lemma}

Lemma \ref{lemma:sms} appears as \cite[7.2]{M06}, while Lemma
\ref{lemma:indiscernibles} is found in \cite[Sec.\ 6]{M06}.

\begin{lemma}
\label{lemma:indiscernibles}
Let $K$ be a differentially closed field and
$$A=\{y\in K: y\ne 0~\&~y\ne 1~\&~y^\prime=y^3-y^2\}.$$
Then $A$ is a strongly minimal set of indiscernibles. 
\end{lemma}

It follows from indiscernibility that $A$ must be a trivial
strongly minimal set and hence $A$ is orthogonal
to each of the sets $E^\sharp_a$.  (Also, the set $A$ is computable in the Turing degree
of the differential field $K$, as defined in the next subsection.)

\begin{lemma}
If $a,b,c,d,\in A$, $a\ne b$, $c\ne d$ and $\{a,b\}\ne \{c, d\}$, 
then $a+b$ and $c+d$ are algebraically independent.
\comment{
ii) If $a,b,c,d\in I$, $a\ne b$, $c\ne d$ and $\alpha$ is differentially transcendental
over $a,b,c,d$ then $a+b$ and $c+d+\alpha$ are algebraically independent.
}
\end{lemma} 
\begin{pf}
Suppose $p(X,Y)\in \Q[X,Y]$ such that $p(a+b,c+d)=0$.  There are only finitely many $y$
with $p(a+b,y)=0$. Suppose without loss of generality that $d\not\in\{a,b\}$.
Then by indiscernibility $p(a+b,c+e)=0$ for every $e\in A\setminus\{a,b,c\}$, a contradiction.
\comment{
ii) If $p(a+b, \alpha+c+d)=0$ then $\alpha$ is algebraic over $a,b,c,d$, a contradiction. 
}
\qed\end{pf}

\subsection{Background in Computable Model Theory}
\label{subsec:CMT}

Now we describe the necessary concepts from computable model theory.
For Proposition \ref{prop:DCFtoGraph} and
Theorem\ref{thm:low}, only Definition \ref{defn:spectrum} is essential,
but the rest of the subsection will make clear why
the broad results in Section \ref{sec:conclusion} are of interest.

Let $\S$ be a first-order structure on the domain $\omega$,
in a computable language (e.g., any language with finitely many
function and relation symbols).  The \emph{(Turing) degree}
$\deg{\S}$ is the Turing degree of the atomic diagram of $\S$;
in a finite language, this is the join of the degrees of the functions
and relations in $\S$.  $\S$ is \emph{computable} if this degree
is the computable degree $\bfz$.  A structure isomorphic to
a computable structure is said to be \emph{computably presentable};
many countable structures fail to be computably presentable.
A more exact measure of the presentability of (the isomorphism
type of) the structure is given by its \emph{Turing degree spectrum}.
\begin{defn}
\label{defn:spectrum}
The \emph{spectrum} of a countable structure $\S$ is the set
of all Turing degrees of \emph{copies} $\M$ of $\S$:
$$ \set{\deg{\M}}{\M\cong\S~\&~\dom{\M}=\omega}.$$
\end{defn}
When dealing with fields, we often write
$\{ x_0,x_1,\ldots\}$ for the domain; otherwise the element
$1$ in $\omega$ might easily be confused with the multiplicative
identity in the field, for instance.  In \cite{K86}, Knight proved that spectra
are always closed upwards, except in a few 
``automorphically trivial'' cases (such as the complete graph on countably
many vertices, whose spectrum is $\{\bfz\}$).

A wide range of theorems is known about the possible spectra
of specific classes of countable structures.  Many classes,
including directed and undirected graphs, partial orders, lattices,
nilpotent groups (see \cite{HKSS02} for these results), and fields
(see \cite{MPSS14}), are known to realize all possible spectra.
We will use the following theorem of Hirschfeldt, Khoussainov,
Shore, and Slinko.
\begin{thm}[see Theorem 1.22 in \cite{HKSS02}]
\label{thm:HKSS}
For every countable, automorphically nontrivial structure $\M$
in any computable language, there exists a (symmetric, irreflexive)
graph with the same spectrum as $\M$.
\end{thm}

Richter showed in \cite{R81} that linear orders, trees and Boolean algebras
fail to realize any spectrum containing a least degree under Turing reducibility,
except when that least degree is $\bfz$, whereas undirected graphs can
realize all such spectra.  Boolean algebras were then
distinguished from these other two classes when Downey and Jockusch
showed that every low Boolean algebra has the degree $\bfz$ in its spectrum;
this has subsequently been extended as far as low$_4$ Boolean algebras,
in \cite{DJ94,KS00,T95}.  In contrast, Jockusch and Soare showed in \cite{JS91}
that each low degree does lie in the spectrum of some linear order with
no computable presentation, although it remains open whether
there is a single linear order whose spectrum contains all
nonzero degrees but not $\bfz$.  (There does exist a graph
whose spectrum contains all degrees except $\bfz$, by results
in \cite{S98,W98}.  A useful survey of related results appears in \cite{FHKKM12}.)

Of relevance to our investigations are the algebraically closed fields,
the models of the closely related theories \ACF0\ and $\textbf{ACF}_p$.
Here the spectrum question has long been settled:
every countable algebraically closed field has every Turing degree
in its spectrum.  On the other hand, every field becomes a constant
differential field when given the zero derivation, which adds no computational
complexity, and so the result from \cite{MPSS14} for fields,
mentioned above, shows that
every possible spectrum is the spectrum of a differential field.
These bounds leave a wide range of possibilities for spectra
of differentially closed fields, and this is the subject of the present paper.
It should be noted that, although every differentially closed field $K$
is also algebraically closed and therefore is isomorphic (as a field) to
a computable field, it may be impossible to add a computable derivation
to the computable field in such a way as to make it isomorphic
(as a differential field) to $K$.

We will show in Proposition \ref{prop:DCFtoGraph} that countable differentially
closed fields do realize a substantial number of quite nontrivial spectra,
derived in a straightforward way from the spectra of undirected graphs.
In particular, differentially closed fields can have all possible proper
$\alpha$-th jump degrees (as defined in that section), for all computable
ordinals $\alpha>0$.  Section \ref{sec:ENI} is devoted to general background
material for the proof of Proposition \ref{prop:DCFtoGraph}.
On the other hand, we then prove Theorem \ref{thm:low},
paralleling the original Downey-Jocksuch result:
it shows that if the spectrum of a countable model of \DCF0 contains
a low degree, then it must also contain the degree $\bfz$.  \DCF0
thus becomes the second theory known to have this property
(apart from trivial examples such as \ACF0).  Our positive results
in Section \ref{sec:spectra}, however, show that this theorem does not
extend to low$_2$ degrees, let alone to low$_4$ degrees, as holds
for Boolean algebras.  Thus \DCF0 realizes a collection of spectra
not currently known to be realized by the models of any other theory
in everyday mathematics.  Finally, in Section \ref{sec:conclusion},
we relativize Theorem \ref{thm:low} and combine it with the results
from Section \ref{sec:spectra} to characterize the spectra of models of \DCF0
precisely as the preimages under the jump operation of the spectra of
automorphically nontrivial graphs, and also as those spectra of such graphs
which are closed under first-jump equivalence.

\section{Eventually Non-isolated Types}
\label{sec:ENI}

The model-theoretic basis of Proposition \ref{prop:DCFtoGraph}
is ENI-DOP, the \emph{Eventually Non-Isolated Dimension
Order Property}, developed by Shelah \cite{SHM84} in proving
Vaught's Conjecture for $\omega$-stable theories.  In this section
we give a simple example of how this property can be used
to code graphs into models of theories satisfying ENI-DOP.
The example may help demystify the coding in
Section \ref{sec:spectra}, which is a more complicated example
of the same phenomenon

In our simple example, we have a language with
two sorts $A$ and $F$, and three unary function symbols
$\pi_1,\pi_2:F\to A$, and $S:F\to F$.
Our theory $T$ includes axioms saying that $A$ is infinite,
the map $(\pi_1,\pi_2):F\to A^2$ is onto,
$\pi_i\circ S=\pi_i$, and $S$ is a permutation of $F$
with no cycles.
This $T$ is complete and has quantifier elimination.
Its prime model consists of a countable set $A$ with one
$\Z$-chain $F_{ab}$ (under $S$) in $F$ for each pair
$(a,b)\in A^2$.  ($F_{ab}$ is the preimage
$(\pi_1,\pi_2)^{-1}(a,b)$ 
and is called the \emph{fiber above} $(a,b)$.)
Every permutation
of $A$ extends to an automorphism of the prime model,
and so $A$ is a set of indiscernibles, in this model
and also in every other model of $T$.

The type over $a$ and $b$ of a single element $x$ of the
fiber $F_{ab}$ is isolated by the formula $(\pi_1(x)=a~\&~\pi_2(x)=b)$.
However, over one realization $c$ of this type, the type
of a new element of $F_{ab}$ (not in the $\Z$-chain of $c$)
over $a$, $b$, and $c$ is not isolated.  This makes the type
of $x$ over $a$ and $b$ an example of an \emph{eventually
non-isolated type}: over sufficiently many realizations of itself,
the generic realization of the type is non-isolated.

The important point here is that
we can add a new point to $F_{ab}$ without forcing any new points
to appear in any other fiber or in $A$.  (Indeed, we can continue adding points
to various fibers without ever forcing any unintended points to appear
in other fibers or in $A$.)  This is roughly what is meant by saying that
the types of generic elements of distinct fibers are \emph{orthogonal}.

We use dimensions to code a graph $G$ on $A$ into a model of this theory $T$.
(The \emph{dimension} of $F_{ab}$ is the number of $\Z$-chains in $F_{ab}$.)
Starting with the prime model of $T$, we add one new $\Z$-chain
to each fiber $F_{ab}$ for which
the graph has an edge between $a$ and $b$.  The orthogonality
ensures the accuracy of this coding, by guaranteeing that this process
does not accidentally give rise to new elements in any fiber $F_{ab}$
for which the graph had no edge between $a$ and $b$.  This builds
a new model $\M$ of $T$, and the permutations of $A$ which extend
to automorphisms of $\M$ are exactly the automorphisms of $G$.

It now follows that there exist continuum-many countable
pairwise non-isomorphic models of $T$, since an isomorphism $f$
between two such structures $\mathfrak{A}$ and $\mathfrak{B}$
would have to map the set of indiscernibles in $\mathfrak{A}$
onto that in $\mathfrak{B}$, hence likewise for the fibers,
and therefore $f$ on the indiscernibles would define an isomorphism
between the graphs coded into $\mathfrak{A}$ and $\mathfrak{B}$.
Moreover, the graph $G$ coded into $\mathfrak{A}$ can be recovered
from the computable infinitary $\Sigma_2$-theory of $\mathfrak{A}$
-- that is, we can compute a copy of $G$ if we know this theory --
and in fact we can enumerate the edges in a copy of $G$ just from the
computable infinitary $\Pi_1$-theory of $\mathfrak{A}$, since this much
information allows us to recognize any two elements of $F_{ab}$ in
$\mathfrak{A}$ that realize the nonisolated $2$-type.

We will use this same strategy to code graphs into countable models
$K$ of \DCF0, using the set $A$ of indiscernibles given by Lemma
\ref{lemma:indiscernibles}.  The fiber $F_{mn}$ for $a_m,a_n\in A$
will be the Manin kernel $E^\#_{a_ma_n}$, defined
in Theorem \ref{thm:Maninkernel} and shown in Theorem
\ref{thm:Maninprops} to have the appropriate properties, and
the non-isolated computable infinitary $\Pi_1$-type
in $F_{mn}$ will be the type of an element of $F_{mn}$
whose coordinates are both transcendental over $\Q\la a_m+a_n\ra$.
With this background, the reader should be ready to proceed
with Proposition \ref{prop:DCFtoGraph}.

Although we will not attempt to generalize here, it is reasonable
to conjecture that the procedure in Section \ref{sec:spectra} should
work for other classes of countable structures for which
similar conditions hold.  Analogues of its converse
(Theorem \ref{thm:low}, essentially) for such classes may be
more challenging.

\comment{
Indeed, if the conditions hold for
types using computable infinitary $\Pi_n$-formulas,
then we conjecture that the same procedure allows one
to code a graph $G$ into a structure $\mathfrak{A}$ in $\C$
in such a way that the computable infinitary $\Pi_n$-theory
of $\mathfrak{A}$ allows one to enumerate the edges in a copy
of $G$.  In this case, the spectrum of $\mathfrak{A}$ ought
to contain exactly those Turing degrees $\bfd$ whose $n$-th
jump $\bfd^{(n)}$ can enumerate the edges in a copy of $G$.
That is, $\spec{\mathfrak{A}}$ should be the preimage of $\spec{H}$
under the $n$-th jump operator (for the $H$ defined from $G$
in Lemma \ref{lemma:GraphEnumSpec} below).
On the other hand, there is no obvious reason why Theorem
\ref{thm:low} need hold for countable models of such a theory.
\DCF0 may be unusual in possessing both ENI-DOP (witnessed
by $\Pi_1$-computable formulas)
and the property that all of its low$_1$ models are computably presentable.
}

\section{Noncomputable Differentially Closed Fields}
\label{sec:spectra}

In this section we consider countable models of the theory \DCF0
which have no computable presentations.  Using countable graphs
with known spectra, we show how to construct differentially closed fields
with spectra derived from those of the graphs.  In particular, we
create numerous countable differentially closed fields which are
not computably presentable.  We show that models of \DCF0 can have
proper $\alpha$-th jump degree for every computable nonzero ordinal
$\alpha$.  However, we will see in Section \ref{sec:conclusion}
that this is impossible when $\alpha=0$:
no countable model of \DCF0 can have a least degree in its
spectrum, unless that degree is $\bfz$.
We encourage the reader to review Section \ref{sec:ENI}
in order to understand the framework for the proof
of the following theorem.

\begin{prop}
\label{prop:DCFtoGraph}
Let $G$ be a countable symmetric irreflexive graph.
Then there exists a countable differentially closed
field $\Khat$ of characteristic $0$ such that 
$$\spec{\Khat}=\set{\bfd}{\bfd'\text{~can enumerate a copy of~}G}.$$
(Saying that a degree $\bfc$ can \emph{enumerate a copy of $G$}
means that there is a graph on $\omega$, isomorphic to $G$, whose
edge relation is $\bfc$-computably enumerable.)
\end{prop}
\begin{pf}
Taking $G$ to have domain $\omega$, we first describe one
presentation of $\Khat$, on the domain $\omega$, without regard to effectiveness.
We begin with $\Qhat$, the differential closure of the constant field $\Q$.
Recall from Subsection \ref{subsec:modeltheory} that
the following is a computable infinite set of indiscernibles:
$$ A=\set{y\in\Qhat}{y'=y^3-y^2~\&~y\neq 0~\&~y\neq 1}.$$
Writing $A=\{a_0<a_1<\cdots\}$,
we use $a_n$ to represent the node $n$ from $G$.

For each $a_m$ and $a_n$ with $m<n$, let $E_{a_ma_n}$ be the elliptic curve
defined by the equation $y^2=x(x-1)(x-a_m-a_n)$.  The type of a differential
transcendental is orthogonal to each strongly minimal set
defined over $\Qhat$.  Thus, for each $m<n$, the Manin kernel
$E^\sharp_{a_ma_n}$ contains only points differentially algebraic
over $\Q\la a_m,a_n\ra $. 
These sets are also orthogonal to $A$.  The points of $E_{a_ma_n}$
in $(\Qhat)^2$ form an abelian group, with (for each $k>0$)
exactly $k^2$ points whose torsion divides $k$, and with
no non-torsion points, since $\Qhat$ is the prime model
of $\textbf{DCF}_0$ over $\Q$.
We will code our graph using these Manin kernels $E^\sharp_{a_ma_n}$,
by adding a new point to $E^\sharp_{a_ma_n}$
(with coordinates transcendental over $\Q\la a_m+a_n\ra$)
to our differential field just if the graph contains an edge from $m$ to $n$.
Any two of these Manin kernels are orthogonal, so adding a point to one 
(or even to infinitely many) of them will not add points to any other.  Similarly,
adding points to the Manin kernels will not add new points to $A$.

\comment{
Let $G$ be an (irreflexive) graph on $I$.
If there is an edge between $a$ and $b$ let $x_{a,b}$ be a new point of $E^\sharp_{a+b}$.
If not let $x_{a,b}$ be a new point of $E^\sharp_{a+b+\alpha}$. 
[Note: the type of a generic point of $E^\sharp_a$ over $k_1$ is given by saying that you are in $E^\sharp_a$ but not algebraic over $\Q\la a\ra $ so easily computable.]
Let $k$ be the differential closure of $k_1\la x_{a,b}: a\ne b\in I\ra $.
Because of the mutual orthogonality of the Manin kernels, $E^\sharp_{a+b}(k)$
contains a point non-algebraic over $\Q\la a+b\ra $ if and only if there is an
edge between $a$ and $b$ if and only if $E^\sharp_{a+b+\alpha}(k)$ contains
only points algebraic over $\Q\la a+b+\alpha\ra $.  Also note that $I(k)=I(k_0)$
so we can easily recover the graph from $k$.
}

Now we build a differential field extension $K$ of $\Qhat$, by adjoining to
$\Qhat$ exactly one new point $x_{mn}$ of $E^\sharp_{a_ma_n}$ for each $m<n$
such that $G$ has an edge between its nodes $m$ and $n$.  
(We note that, by orthogonality, the type of each generic point
of $E^\sharp_{a_ma_n}$ over the differential field $L$ generated
by the preceding points $x_{m'n'}$ is computable:  it is given by saying that $x_{mn}$
is in $E^\sharp_{a_ma_n}$ but is not algebraic over $L\la a_m+a_n\ra $.)
Adjoining all these $x_{mn}$ yields a 
differential field $K$, and the differential field we want is the differential
closure $\Khat$ of this $K$.
The principal relevant feature of $\Khat$
is that, 
because of the mutual orthogonality of the Manin kernels, $\Khat(E^\sharp_{a_ma_n})$
contains a point non-algebraic over $\Q\la a_m+a_n\ra $ if \emph{and only if}
there is an edge between $m$ and $n$ in $G$.

Now we claim that the spectrum of this $\Khat$ contains exactly those Turing degrees
whose jumps can enumerate a copy of $G$.  To show that every degree in the spectrum
has this property, suppose that $L\cong\Khat$ has degree $\bfd$.
Then with a $\bfd$-oracle, we can decide the set of all nontrivial solutions $b_0 < b_1<\cdots$
in $L$ to $y'=y^3-y^2$.  (The only trivial solutions are $0$ and $1$.)
We build a graph $H$, with domain $\omega$,
using a $\bfd'$-oracle.  The oracle tells us, for each $m<n$ and each solution
$(x,y)\in L(E^{\#}_{b_mb_n})$, whether or not $x$ is algebraic over $\Q\la b_m+b_n\ra$.
If so, then we go on to the next point in $L(E^{\#}_{b_mb_n})$.  If $x$ is not algebraic,
then we enumerate an edge between $m$ and $n$ into our graph $H$.
The graph $H$ thus enumerated is isomorphic to $G$:  the isomorphism $f$
from $L$ onto $\Khat$ must map the set $\{b_0,b_1,\ldots\}$ bijectively
onto the set $\{ a_0,a_1,\ldots\}$, and the map sending each $m\in H$
to the unique $n\in G$ with $f(b_m)=a_n$ will be an isomorphism
of graphs.  Thus $\bfd'$ has enumerated a copy $H$ of $G$.

Conversely, suppose that the Turing degree $\bfd'$ enumerates
a graph $H$ isomorphic to $G$.  Specifically, for a fixed set $D\in\bfd$,
there is a Turing functional $\Phi$ for which 
the edge relation on $H^2$
is the domain of the partial function $\Phi^{D'}$.
The description of $\Khat$ above explains how to build a differentially
closed field $\Lhat$ below a $\bfd$-oracle with $\Lhat\cong\Khat$.
Using Theorem \ref{thm:Harrington}, start building a computable
copy of $\Qhat$, in which we enumerate all nontrivial
solutions $b_n$ to $y'=y^3-y^2$, but build this solution slowly,
with one new element at each stage, so that each step $L_s$ in this
construction is actually a finite fragment of the differential
field $L$ we wish to build.  Then, with the $\bfd$-oracle,
enumerate the jump $D'$ of the set $D\in\bfd$:  say $D'=\cup_{s\in\omega}D'_s$.
Whenever we find a stage $s$ such that some $\la m,n\ra$ lies in
$\dom{\Phi^{D'_s}_s}$ (and did not lie in this domain for $s-1$),
we adjoin to $L_s$ a new point $(x_{m,n,s},y_{m,n,s})$ in $E^\sharp_{b_mb_n}$,
such that $x_{m,n,s}$ does not yet satisfy any nonzero differential polynomial
at all over $L_s$, and is specified not to be a zero of the first $s$ polynomials
of degree $\leq s$ over $L_s$.  Of course,
$y_{m,n,s}$ is a zero of the curve $E^\sharp_{b_mb_n}$ over $x_{m,n,s}$; this fully
determines $y_{m,n,s}$ and its derivatives in terms of $L_s$ and $x_{m,n,s}$
and its derivatives.

At the next stage,
if we still have $\la m,n\ra\in\dom{\Phi^{D'_{s+1}}_{s+1}}$,
we declare that $x_{m,n,s+1}=x_{m,n,s}$ is not a zero of any
of the first $s+1$ polynomials of degree $\leq s+1$ over $L_{s+1}$.
If we ever reach a stage
$t>s$ with $\la m,n\ra\notin\dom{\Phi^{D'_{t}}_{t}}$ (which is possible,
if the oracle has changed from the previous stage), then we turn $(x_{m,n,s},y_{m,n,s})$
into a $k$-torsion point, with $k\geq t$ being the smallest value for which this
is consistent with the finite fragment $L_{t-1}$ built up till then.  Since 
the types of torsion points are dense in the space of all types,
the finitely many facts we have enumerated so far
about $L_{t-1}$ cannot possibly force this point to be a non-torsion point,
so for some $k$ this will be possible, and by searching we can identify such a $k$,
using the decidability of \DCF0.
As we subsequently continue to build $L$ (including the cofinite portion
of $\Qhat$ which is yet to be constructed), we will take this $k$-torsion
point into account, treating it as part of $\Qhat$.  The decidability of
\DCF0 makes it easy to include the point into $\Qhat$
and still know what to build at each subsequent step.

Thus the existence of
a nonalgebraic point on $E^\sharp_{b_mb_n}$ in the field $L$ built by this process
is equivalent to $\la m,n\ra$ actually lying in $\dom{\Phi^{D'}}$,
and for all $\la m,n\ra$ not in this domain, every pair $(x_{m,n,s},y_{m,n,s})$
ever defined (for any $s$) was eventually turned into a torsion point,
so that it wound up in the subfield $\Qhat$ of $L$, since this subfield
contains all $k^2$ of the $k$-torsion points for $E_{b_mb_n}$ in $L$.
Therefore, the $L$ that we finally built is the differential field extension
of $\Qhat$ by one nontorsion point for each edge in $H$, hence is
isomorphic to the differential field $K$ built above.  So the differential
closure $\Lhat$ of $L$ is isomorphic to $\Khat$, and
is also $\bfd$-computable, by Theorem \ref{thm:Harrington}.
\qed\end{pf}

We note that in Proposition \ref{prop:DCFtoGraph}, it is reasonable
to replace the graph $G$, which the $\bfd'$-oracle can enumerate,
by another countable graph $H$ which the same oracle can actually compute.
The converse is accomplished by the technique known as a
Marker $\exists$-extension.  The forwards direction too is
a simple question of coding.  

\begin{lemma}
\label{lemma:GraphEnumSpec}
Let $H$ be a countable (symmetric irreflexive) graph.
Then there exists a countable graph $G$ such that
$$\spec{H}=\set{\bfd}{\bfd\text{~can enumerate a copy of~}G}.$$

Conversely, for every countable graph $G$, there exists a countable
graph $H$ satisfying this same equation.
\qed\end{lemma}

Recall that, for a computable ordinal $\alpha$, the $\alpha$-th \emph{jump degree}
of a countable structure $\M$ is the least degree in the set $\set{\bfd^{(\alpha)}}{\bfd\in\spec{\M}}$.

\begin{thm}
\label{thm:jumpspectra}
For every graph $H$, there exists a differentially closed field $K$
such that
$$\spec{K}=\set{\bfd}{\bfd'\in\spec{H}}.$$
In particular, for every computable ordinal $\alpha>0$ and every
degree $\bfc >_T \bfz^{(\alpha)}$, there is a differentially closed
field which has $\alpha$-th jump degree $\bfc$, but has no
$\gamma$-th jump degree whenever $\gamma<\alpha$.
\end{thm}
Using ordinal addition, one can re-express the second result by stating that,
for every $\beta<\omega_1^{CK}$ and every $\bfc$ with $\bfc>_T\bfz^{(1+\beta)}$,
there is a differentially closed field $K$ with proper $(1+\beta)$-th jump degree $\bfc$.

\begin{pf}
Given $H$, use Lemma \ref{lemma:GraphEnumSpec} to get a graph
$G$ whose copies are enumerable by precisely the Turing degrees in $\spec{H}$.
Then apply Proposition \ref{prop:DCFtoGraph} to this $G$
to get the differentially closed field $K$ required, with
$$ \spec{K}=\set{\bfd}{\bfd'\text{~can enumerate a copy of~}G}
=\set{\bfd}{\bfd'\in\spec{H}}.$$

Now, for every computable ordinal $\beta$ and every degree $\bfc\geq \bfz^{(\beta)}$,
there exists a graph $H$ with $\beta$-th jump degree $\bfc$,
but with no $\gamma$-th jump degree for any $\gamma<\beta$.
(This is shown for linear orders in \cite{AJK90} and \cite{DK92}
for all $\beta\geq 2$, and Theorem \ref{thm:HKSS} then transfers the result to graphs.
For $\beta<2$ it is a standard fact; see e.g.\ \cite{FHKKM12}.)
If $\alpha>0$ is finite, let $\beta$ be its predecessor and
apply the first part of the corollary to the $H$ corresponding to
$\bfc$ and to this $\beta$.
Then
$$\set{\bfd^{(\beta)}}{\bfd\in\spec{H}} = 
\set{(\bfd')^{(\beta)}}{\bfd\in\spec{K}} =
\set{\bfd^{(\alpha)}}{\bfd\in\spec{K}},$$
so $\bfc$ is the $\alpha$-th jump degree of $K$.
When $\beta\geq\omega$, the degree
$(\bfd')^{(\beta)}$ is just $\bfd^{(\beta)}$ itself, and so, for every infinite computable
ordinal $\alpha$, the above analysis with $\beta=\alpha$ shows that
again $K$ has $\alpha$-th jump degree $\bfc$.
In both cases, 
this also proves that for each $\gamma<\alpha$,
$K$ has no $\gamma$-th jump degree.
\qed\end{pf}

\section{Low Differentially Closed Fields}
\label{sec:low}

Theorem \ref{thm:jumpspectra} demonstrated that, for every nonlow
Turing degree $\bfd$, there exists a $\bfd$-computable differentially
closed field with no computable presentation:  with $\bfd'>\bfz'$, just take
the model of \DCF0 given by the corollary with jump degree $\bfd'$.
(The corollary showed specifically that every degree whose jump computes $\bfd'$
lies in the spectrum, so the structure has a $\bfd$-computable copy.)
Of course, there exist noncomputable low Turing degrees $\bfd$,
i.e., degrees with $\bfd>\bfz$ but $\bfd'=\bfz'$.  Theorem
\ref{thm:jumpspectra} yields no proof of the same result
for these degrees.  Indeed, the surprising answer is that
when $\bfd$ is low, every $\bfd$-computable
differentially closed field has the degree $\bfz$ in its spectrum.

\begin{thm}
\label{thm:low}
Every low differentially closed field $K$ of characteristic $0$
is isomorphic to a computable differential field.
\end{thm}

Before beginning the full proof, we give some idea how it will go.
Our goal is to construct a computable differential field $F$,
with elements $y_0,y_1,\ldots$, isomorphic to $K$, whose
elements are $x_0,x_1,\ldots$.  The isomorphism $x_n\mapsto y_{h(n)}$
will be $\Delta^0_2$, and we construct finite approximations $h_s$
to $h$.  We must ensure that the limit of these $h_s$
exists and is a bijection.
The requirement $\R_n$ is that $\lim_sh_s(n)$ exists;
the requirement $\S_m$ is that $\lim_s h_s^{-1}(m)$ exists.
Since each $h_s$ will define a finite partial isomorphism
into $F$ from the current approximation $K_s$ to $K$,
the limit $h$ will then define an isomorphism from $K$ onto $F$.

For a single element $x_n\in K$, the basic module for satisfying
$\R_n$ is not difficult.  Since $K$ is low,
we can guess effectively at the minimal differential polynomial $p_n$ of $x_n$
over the finitely many higher-priority elements $x_i$ of $K$.  Assuming at stage
$s$ that our current guess is correct, we simply check through the finitely
many elements currently in $F$ to see whether this $p_n$ is currently the minimal
polynomial of any of them over the corresponding $y_{h_s(i)}$ in $F$.
If so (and if that element is not already claimed by a higher-priority requirement),
then we choose it as the image of $x_n$.  If not, then we add a new element
to $F$, making it a solution of $p_n$, and define it to be the image
of $x_n$ at this stage.  Once our guesses at $p_n$ have stabilized,
this element will be $y_{h(n)}$, the image of $x_n$ under our
$\Delta^0_2$-isomorphism.  In the meantime, if our guess at $p_n$
changes, we simply start the process over, leaving a leftover
element in $F$.

Our construction will define the atomic diagram of $F$ only in ways
consistent with the complete decidable theory \DCF0.  (If $\R_n$
wants $p_n$ to be given a zero in $F$, but \DCF0 refuses to allow it,
then the construction waits for a change in the guesses $p_0,\ldots,p_n$,
which must happen, since $K$ satisfies \DCF0.)  Therefore, a leftover element
still can rely on \DCF0's assurance that there exists some zero of that $p_n$:
$K$ must contain some such zero, although $x_n$ turned out not to be such a zero.
Our next task, in building $F$, is to find a preimage for each leftover element
$y_m$, as required by $\S_m$.  Of course, once $y_m$
is made into a root of a certain differential polynomial, it must
remain a root of that polynomial; however, it might later be made into
a root of another differential polynomial of lesser rank, so that the
first one might not be its \emph{minimal} differential polynomial.
Since ranks are ordinals, this can only happen finitely often.

While $y_m$ is believed to have minimal differential polynomial $f$
over the higher-priority elements of $F$, and while $h_s^{-1}(m)$
is undefined, we search for an element of $K$ which appears
to have the same minimal differential polynomial over the corresponding
higher-priority elements of $K$.  If we find one, we make it the preimage of $y_m$.
However, the existence of such an element in $K$ is guaranteed only
if $f$ is constrainable (over the differential subfield generated by the higher-priority
elements), which may not be decidable.  Moreover, even if we find an $x\in K$
which appears to have the correct minimal differential polynomial,
we could turn out to be mistaken, since we have only a computable approximation
to minimal differential polynomials in $K$.  There is a danger that no
$x$ with the correct minimal differential polynomial actually exists in $K$,
but that $K$ keeps offering us different possible elements $x$ forever, each appearing
to have the minimal differential polynomial we want.  (In this sense, $K$
``cannot be trusted'' ever to give us a correct preimage, nor to cease supplying
possibilities which turn out to be incorrect.)  Therefore, while searching
for a zero of $f$ in $K$, we use the type function $\gamma$ from
Lemma \ref{lemma:typefct} to determine a principal type containing
the formula $f=0$, and make this the type of $y_m$.
The ground field (providing the $e$ in Lemma \ref{lemma:typefct})
is the differential subfield of $F$ generated by
the higher-priority elements, under the assumption that no higher-priority
requirement ever acts again.  Obeying the type function ensures
that eventually $y_m$ will settle as a zero of a polynomial
which is constrainable (over the higher-priority elements of $F$),
and this in turn ensures that $K$ will contain an element with
that same minimal differential polynomial, which we will eventually
find and define to be the preimage of $y_m$.
The construction is therefore a finite-injury procedure,
using these basic modules for the two types of requirements.

\begin{pf}
Our goal is to build a computable differential field $F$,
with domain $\{ y_0,y_1,\ldots\}$, and a sequence of uniformly
computable finite partial functions $h_s:\omega\to\omega$ such that,
for all $n$, $\lim_s h_s(n)$ converges to an element $h(n)$
so as to define an isomorphism $x_n\mapsto y_{h(n)}$ from $K$ onto $F$.
When $n\leq h(n)$, we will arrange that $x_n$ and $y_{h(n)}$
have the same minimal differential polynomials over the  differential subfields
generated by the higher-priority elements in $K$ and $F$:
\begin{align*}
\Q\la x_0, x_{h^{-1}(0)},x_1,x_{h^{-1}(1)},&\ldots,x_{n-1},x_{h^{-1}(n-1)}\ra\subseteq K\\
\Q\la y_{h(0)},~y_0,~y_{h(1)},~y_1,~&\ldots,y_{h(n-1)},~y_{n-1}\ra\subseteq F.
\end{align*}
More precisely, there will be a differential polynomial $p_n\in\Q\{ X_0,Y_0,X_1,\ldots,Y_{n-1},X_n\}$
such that $p_n(x_0, x_{h^{-1}(0)},x_1,x_{h^{-1}(1)},\ldots,x_{h^{-1}(n-1)},X_n)$
is the minimal differential polynomial of $x_n$ over the first subfield
and $p_n(y_{h(0)},y_0,y_{h(1)},y_1,\ldots,y_{n-1},Y_n)$ is the minimal
differential polynomial of $y_{h(n)}$ over the second subfield.

Likewise, when $n > h(n)$, we will arrange that $x_n$ and $y_{h(n)}$ have
the same minimal differential polynomials over the differential subfields
generated by higher-priority elements:
\begin{align*}
\Q\la x_0, x_{h^{-1}(0)},x_1,x_{h^{-1}(1)},&\ldots,x_{h^{-1}(h(n)-1)},x_{h(n)}\ra\subseteq K\\
\Q\la y_{h(0)},~y_0,~y_{h(1)},~y_1,~&\ldots,y_{h(n)-1},~y_{h(h(n))}\ra\subseteq F.
\end{align*}
(With $n>h(n)$, the lower index $h(n)$ gives the priority of the pair $(x_n,y_{h(n)})$.
Those pairs containing any of the elements $x_0,\ldots,x_{h(n)}$ and
$y_0,\ldots,y_{h(n)-1}$ will have higher priority and so will be considered first.)

This will establish that $h$ defines an embedding of differential fields.
We will also ensure that $h:\omega\to\omega$ is a bijection,
hence defines an isomorphism.  Since $F$ is computable, this will
prove the theorem.  
\comment{
Notice that, in contrast to the situation with Boolean algebras,
it will follow that every low differentially closed field is $\Delta^0_2$-isomorphic
to a computable one; for Boolean algebras a $\Delta^0_3$-isomorphism is
sometimes required.
}

Our key asset in this construction is a computable approximation
not only of the atomic diagram of the differential field $K$, but also of the minimal
differential polynomial of each element $x_n$ (in the domain
$\{ x_0,x_1,\ldots\}$ of $K$) over the differential subfield
$\Q\la x_0,\ldots,x_{n-1}\ra$.  Indeed, we use slightly more:
we can effectively approximate the minimal differential
polynomial $p_{n,\rho}$ of any $x_n$ over $\Q\la x_{n_1},\ldots,x_{n_k}\ra$,
where $\rho=\la n_1,\ldots,n_k\ra\in\omega^{<\omega}$.
This holds because the computable infinitary
$\Sigma_1$-diagram of $K$ is computable in the jump $(\deg{K})'$,
i.e., in $\bfz'$.  Recall that by our convention in this article,
the minimal differential polynomial of a differential transcendental is
the zero polynomial, 
and the comments above apply to differential transcendentals as well,
since one jump over $\deg{K}$ is enough to decide whether $x_n$
satisfies any nonzero differential polynomial at all over $\Q\la x_{n_1},\ldots,x_{n_k}\ra$.

So, for each $n$,
we will guess at some $p_n\in\Q\{ X_0,Y_0,\ldots,Y_{n-1},X_n\}$ giving
the minimal differential polynomials of $x_n$ and $y_{h(n)}$ over the
relevant differential subfields, as described earlier.
Our requirements to satisfy are:
\begin{align*}
\R_n: & ~h(n)=\lim_s h_s(n)\text{~exists}\\
\S_m: & ~h^{-1}(m)=\lim_s h_s^{-1}(m)\text{~exists,}
\end{align*}
with priority $\R_0\prec\S_0\prec\R_1\prec\cdots$.
If we can satisfy them, and maintain our rule that each $p_n$
gives the minimal differential polynomial of both $x_n$ and $y_{h(n)}$,
then we will have built our isomorphism, which will in fact then be
$\bfz'$-computable itself.

\comment{
At each step in building $F$, we will ensure that the finite portion
of its atomic diagram so far defined is consistent with the (complete, decidable)
theory $\DCF0$:  if $\Delta_s$ is the portion defined by stage $s$,
using the elements $y_0,\ldots,y_m$, and we wish to adjoin a new
atomic formula $\delta(y_0,\ldots,y_{m+1})$, we simply ask whether
the finitary existential formula
$$ (\exists Y_0\cdots\exists Y_{m+1})[\Delta_s~\&~\delta~\&~
(\forall i\leq m) Y_i\neq Y_{m+1}]$$
lies in $\DCF0$.  If not, then we refuse to add any new element
$y_{m+1}$ satisfying $\delta$ to $F$; in this case the computable
approximation to $K$ which suggested that $K$ satisfied this
existential formula must subsequently change, so we simply
keep approximating $K$ and wait for the change.
Notice that
all formulas here are finitary.  $\DCF0$ cannot decide the consistency
of computable infinitary formulas, so cannot be used to decide
whether a pair of differential polynomials is a constrained pair.
}

The strategy for satisfying a single requirement $\R_n$ is relatively
simple.  There exists a stage $s$ by which our approximation to $K$
will have settled on the
true minimal differential polynomial $p_n(x_0, x_{h_s^{-1}(0)},\ldots,
x_{h_s^{-1}(n-1)},X_n)$ of $x_n$ over the higher-priority elements.
If there already exists an element $y_m$ in $F$ for which
$p_n(y_{h_s(0)},y_0,\ldots,y_{n-1},Y_n)$ is the minimal differential polynomial
(over these higher-priority elements, according to the structure of $F$ at this
stage), we define $h_{s+1}(n)=m$.  (This includes the situation where $m<n$ and
$h_s(n)=m$ was already defined for the sake of the higher-priority $\S_m$.)
Alternatively, if for some lower-priority $y_m$ already in $F$ it is consistent with $\DCF0$
(given the current types of higher-priority elements of $F$), for $y_m$
to become a zero of this polynomial, then again we define $h_{s+1}(n)=m$;
otherwise, we add a new element $y_m$ to $F_{s+1}$, making it a zero
of this polynomial (provided this is consistent, the same as above) and set
$h_{s+1}(n)$ equal to this new $m$.  (If neither of these options
is consistent, then we simply wait for our approximations to $K$ to change.)
Assuming that no higher-priority
requirement ever again injures $\R_n$, and that the guess $p_n$ never
again changes, this $y_m$ will continue to have this
minimal differential polynomial throughout the rest of the construction:
neither $\R_n$ nor any higher-priority requirement will ever need to change it,
and no lower-priority requirement will ever be allowed to do so.
($K$ itself witnesses that it is consistent with \DCF0 for $x_0,x_{h^{-1}(0)},\ldots,x_n$
to have the minimal differential polynomials that we have found,
so \DCF0 will not require any further changes to $y_m$.)
Therefore $\R_n$ will never again injure any lower-priority requirement.
Also, any similar actions taken by $\R_n$ before we reached this stage $s$
will not impede us from satisfying $\R_n$ or any higher-priority requirement.
The strategy for satisfying a single requirement $\S_m$ is more
complicated; we will describe it in the construction, before the instructions
for the $\S$-substages.

\begin{notation}
\label{notation:stages}
To avoid cumbersome subscripts, we adopt the convention of writing
``~\!$[s]$'' at the end of an expression to indicate that all items in the
expression have the values assigned to them as of stage $s$.
For example, $p_{n_i,\rho_i}(y_{h(n_0)},\ldots,y_{h(n_i)})[s]$ will 
denote $p_{n_{i,s},\rho_{i,s},s}(y_{h_s(n_{0,s})},\ldots,y_{h_s(n_{i,s})})$.
\end{notation}

Having $F_s$ be a finite fragment
of a differential field will allow us to lean heavily on the theory \DCF0 for guidance
in constructing $F_{s+1}$.  This theory is complete and decidable, and so, given
the finite fragment $F_s$ containing (say) $y_0,\ldots,y_r$,
we can write out the entire relational atomic diagram $\psi(y_0,\ldots,y_r)$
of these elements.  When considering how
to build $F_{s+1}$, we can then ask whether \DCF0 contains the sentence
$$ \exists Y_0\cdots\exists Y_m [\psi(Y_0,\ldots,Y_r)~\&~g(Y_0,\ldots,Y_m)=0].$$
(Here $g$ is some polynomial over $\Q$
for which we might wish to declare $\yvec$ to be a zero.)  If this is
inconsistent, then the decision procedure for \DCF0 will tell us so,
and we will not set $g(\yvec)=0$ in $F_{s+1}$.  If it is consistent,
it belongs to the complete theory \DCF0, so some tuple
of elements of $K$ must realize $[\psi(\Xvec)~\&~g(\Xvec)=0]$,
and it is safe to set $g(\yvec)=0$ in $F_{s+1}$, as $K$ must contain
preimages of these elements which are consistent with the minimal
differential polynomials $p_0,p_{h^{-1}(0)},p_1,\ldots$, up to the
first $p_n$ for which our approximations have not yet converged.  (Notice that
all formulas here are finitary.  $\DCF0$ cannot decide the consistency
of computable infinitary formulas, so cannot be used to decide, for instance,
whether a pair of differential polynomials is a constrained pair.)  Of course,
we must also verify that doing so will not change the minimal differential
polynomial of any higher-priority element.  Part of the purpose of Lemma
\ref{lemma:minpoly} is to show how to do this verification effectively.

\comment{

The requirements for the construction concern the function $h:\omega\to\omega$
which we build as the limit $\lim_s h_s$ of a computable sequence of functions.
In order to make $h$ a bijection, we will satisfy for each $m$ and $n$:
\begin{align*}
\R_n: & ~h(n)=\lim_s h_s(n)\text{~exists.}\\
\S_m: & ~h^{-1}(m)=\lim_s h_s^{-1}(m)\text{~exists.}
\end{align*}
Additionally, to make $h$ define an isomorphism of differential fields, we will 
build a sequence $\la n_i\ra_{i\in\omega}$, where $i\mapsto n_i$
is a bijection from $\omega$ onto itself, and ensure,
for each $i$ and each $p\in\Q\{ X_0,\ldots,X_i\}$:
$$p(x_{n_0},\ldots,x_{n_i})=0\text{~in~}K\iff p(y_{h(n_0)},\ldots,y_{h(n_i)})=0\text{~in~}F.$$
It will then follow, by induction on $i$, that each $x_{n_i}$ satisfies the same
$1$-type over $K_{i-1}=\Q\la x_{n_0},\ldots,x_{n_{i-1}}\ra$ that $y_{h(n_i)}$ satisfies over
$\Q\la y_{h(n_0)},\ldots,y_{h(n_{i-1})}\ra$, so that $h$ will define an isomorphism.
These requirements are given a priority ranking, with
$\R_i\prec \S_i \prec\R_{i+1}$ for all $i$ (meaning that among these three,
$\R_i$ has the highest priority, then $\S_i$, then $\R_{i+1}$).

}

At stage $0$, we set $F_0$ to contain $y_0=0$ and $y_1=1$
as the identity elements of $F$.
\comment{
It would be natural to define triples such as $(y_0,y_0,y_0)$
and $(y_0,y_1,y_1)$ to lie in the addition relation
for $F_0$, and $(y_0,y_0)$ and $(y_1,y_0)$ to lie in the differentiation relation, and so on.
There would be no harm in doing so for finitely many tuples,
but, in line with the construction of the rest of $F$,
}
The actual step is that we add $Y_0$ and $(Y_1-1)$ to the set $U_0$,
i.e., to the computable enumeration of the set $U$ of those differential polynomials
$f\in\Q\{ Y_0,Y_1,\ldots\}$ for which $f(y_0,y_1,\ldots)=0$ in $F$.
This is equally strong and will simplify the construction, since it parallels our process
for approximating $K$, which uses minimal differential polynomials
rather than using the relations directly.
In order to use the differential polynomials this way, we will
need to be able to consider the finite set $U_s$ at each stage and decide,
for each $m$, just what minimal polynomial (over the higher-priority
elements of $F$) we have committed $y_m$ to satisfy.  This requires the following lemma.

\begin{lemma}
\label{lemma:minpoly}
There is an algorithm which, when given as input (strong indices for) finite sets
$V,W\subseteq\Q\{ T_0,\ldots,T_r\}$ of differential polynomials and an $m\leq r$
such that $\exists T_0\cdots \exists T_r\psi$ lies in \DCF0,
where $\psi$ is the formula
$$ \bigwedge_{g\in V}g(T_0,\ldots,T_r)=0~\&~
\bigwedge_{g\in W}g(T_0,\ldots,T_r)\neq 0~\&~
\bigwedge_{i<j\leq r}T_i\neq T_j,$$
outputs a differential polynomial $f=\sum_\theta f_\theta T_m^{\theta}$
in $\Q\{T_0,\ldots,T_m\}$ of least possible rank in $T_m$
(written here using finitely many $f_\theta\in\Q\{ T_0,\ldots,T_{m-1}\}$)
such that \DCF0 contains the sentence
$$\left(\forall T_0,\ldots\forall T_r [\psi\to~f=0]\right)~\&~
\left( \exists T_0,\ldots\exists T_r \bigvee_\theta [\psi~\&~f_\theta\neq 0]\right).$$
\end{lemma}
(The point here is that committing ourselves to the finite set $\psi$
of conditions will force $T_m$ to be a zero of $f$, but will not force
it to be a zero of any differential polynomial of lesser rank.  So the
algorithm is producing the apparent minimal differential polynomial $f$
of $T_m$ over $T_0,\ldots,T_{m-1}$, under the condition $\psi$,
although of course $\psi$ does not necessarily rule out the possibility
of $T_m$ satisfying some differential polynomial of smaller rank as well.
The $f$ produced is unique up to a scalar from $\Q^{\times}$.)
\begin{pf}
This is simply the algorithm originally developed by Ritt for reducing
one differential polynomial modulo others of lower rank.  It is given in full
in \cite{R32}, in a version which allows for several derivations, and
is analogous to the reduction procedure
for finding a principal generator of an algebraic ideal in the (non-differential)
polynomial ring $L[T]$.  Here we first convert the negative statements
given by $W$ to positive ones by adjoining variables $S_g$ satisfying
$1=S_g\cdot g(T_0,\ldots,T_r)$ for each $g\in W$.  Then we do Ritt's procedure,
using all polynomials in $V$ and these new equations from $W$,
to get a minimal differential polynomial for $T_0$.  If this polynomial
lies in $\Q\{ T_0\}$, then it is our output $f_0$ for the $m=0$ case; if not,
then $f_0$ is the zero polynomial.  In either case, we then treat
the quotient field of $\Q\{ T_0\}/\{ f_0\}$ as our ground field and repeat the
process for $T_1$ over this ground field (still using all the equations
from $V$ and $W$) to produce $f_1$, then continue recursively up to $f_m$
which is the desired $f$.
\qed\end{pf}

Now we give the algorithm to be followed at stage $s+1$,
using the function $h_s$ and the set $U_s$ from stage $s$.
The domain of $h_s$ contains finitely many elements of $\omega$,
which we view as indices of the elements $x_n$ of $K$,
while its range is a set of certain indices $m\leq r$ of
elements $y_m$ of the finite set $F_s=\{ y_0,\ldots,y_r\}$.
We order the indices of elements of $F_s$ according to priority:
$$ h_s(0) \prec 0 \prec h_s(1) \prec 1 \prec\cdots\prec r,$$
and, after removing all repetitions from this list, we name these
indices $m_{0,s}\prec m_{1,s}\prec\cdots$.  
\comment{
That is, $m_{0,s}=h_s(0)$
is the highest-priority index, since keeping it fixed at subsequent stages
will satisfy $\R_0$.  The next highest-priority index is usually $0$ itself,
since $\S_0$ wants $h_s^{-1}(0)$ to stay fixed at all subsequent stages;
however, if $h_s(0)=0$, then $m_{0,s}=0$ and we do not repeat $0$
in our list, but instead set $m_{1,s}=h_s(1)$, since $\R_1$ has the
next highest priority.  
}
If $h_s(n)$ is undefined for some $n$,
we simply skip that spot in our list of indices $m_{i,s}$.
The list ends once it contains all indices of elements of $F_s$,
namely $\{ 0,1,\ldots,r\}$.
For each $i$, we define $n_{i,s}=h_s^{-1}(m_{i,s})$, if this inverse image
exists.  For the least $j$ such that $n_{j,s}$ is not defined by this process,
we set $n_{j,s}$ to be the least element not in $\dom{h_s}$,
since we might be able to extend $\dom{h_{s+1}}$ to include this element.
Then, for each $i\leq j$, we set $\rho_{i,s}$ to be the finite tuple
$(n_{0,s},n_{1,s},\ldots,n_{i-1,s})$ containing those elements
of higher priority than $n_{i,s}$ in $F_s$.

The atomic diagram of $F_s=\{ y_0,\ldots,y_r\}$ so far determined is denoted
$$\psi_{s}(Y_0,\ldots,Y_r):~\bigwedge_{i<j\leq r}
Y_i\neq Y_j~\&~\bigwedge_{f\in U_s}f(Y_0,\ldots,Y_k)=0~\&~
\bigwedge_{i<s\&g_i\notin U_s}g_i(Y_0,\ldots,Y_k)\neq 0.$$
(At the Final Step of each stage $s+1$, it is determined whether
the $s$-th polynomial $g_s$ lies in $U$ or not.)
Similarly, for each $i$ with $n_{i,s}$ defined, $\sigma_{i,s}$ is the current
approximation to $K$ up to $x_{n_{i,s}}$, using the priority ordering:
$$ \sigma_{i}(X_{n_{0}},\ldots,X_{n_{i}})[s]:
\bigwedge_{j\leq i} \left[p_{n_{j},\rho_{j}}(X_{n_{0}},\ldots,
X_{n_{j}})=0~\&~\bigwedge_{k<j}X_{n_{k}}\neq X_{n_{j}}\right] [s]$$
where, as defined earlier, $p_{n_{i},\rho_{i}}(X_{n_{0}},\ldots,X_{n_{i}})[s]$
is the current approximation
to the minimal differential polynomial of $x_{n_{i,s}}$
over $\Q\la \rho_{i,s}\ra$.  (Having $p_{n,\rho,s}$ be the
zero polynomial when $x_n$ appears to be differentially transcendental
over $\Q\la\rho\ra$ suits this definition of $\sigma_{i,s}$ nicely.)

\textbf{$\R_n$-substages.}
At stage $s+1$, we go through each $\R_n$ and $\S_n$ with $n\leq s$ in turn,
with one substage for each, starting with $\R_0$.  At the substage for a requirement $\R_n$,
fix $i$ such that $n=n_{i,s}$.  (Such an $i$ must exist,
since we included the least index $\notin\dom{h_s}$
on our list of indices $n_{i,s}$.  After this least
element has been reached, no further substages will be executed
at this stage.)
Now we know that, for all $n_{k,s}$ with $k<i$,
$h_{s+1}(n_{k,s})=h_s(n_{k,s})$, since otherwise the stage would
have ended already.  First we check whether the sentence
$$\exists X_{n_{0}}\cdots\exists X_{n_{i}}~\sigma_{i}(X_{n_{0}},\ldots,X_{n_{i}})[s]$$
belongs to \DCF0.  If not, then we do nothing at this substage,
and do not go on to the next substage, but instead go straight
to the Final Step of stage $s+1$ (described below).
In particular, $h_{s+1}(n_{k,s})$ is undefined for all $k\geq i$.
As a simple example, if $p_{n_{i},\rho_{i}}=X_{n_{i}}-a[s]$
and $p_{n_{j},\rho_{j}}=X_{n_{j}}-a[s]$ for the same rational $a$
and for some $j<i$, then the sentence would be rejected as inconsistent.
If it is consistent, then we follow these instructions.
\begin{enumerate}

\item
If $h_{s+1}(n_{i,s})$ has been defined at an earlier substage,
then we keep that value and go on to the next substage.
(This happens if $h_{s+1}(n_{i,s})<n_{i,s}$.)

\item
If $h_s(n_{i,s})\converges$ and Lemma \ref{lemma:minpoly} shows
the minimal differential polynomial of $y_{h_s(n)}$
over $\Q\la y_{h(n_{0})},\ldots,y_{h(n_{i-1})}\ra[s]$
to be $p_{n,\rho_{i}}(y_{h(n_{0})},\ldots,y_{h(n_{i-1})},X)[s]$,
then we preserve the map, setting $h_{s+1}(n)=h_s(n)$, and go on to the
next substage.  For instance, we do this if $\rho_{i,s-1}=\rho_{i,s}$
and $p_{n,\rho_{i}}[s-1]=p_{n,\rho_{i}}[s]$.

\item
Otherwise, either $h_s(n)$ is undefined, or else
$h_s(n)=m'$ is defined with $m'\geq n$ but
$p_{n,\rho_{i}}(y_{h(n_0)},\ldots,y_{h(n_{i-1})},X)[s]$ is not the minimal differential
polynomial of $y_{m'}$ over $\Q\la y_{h(n_0)},\ldots,y_{h(n_{i-1})}\ra[s]$ in $F_s$.
(This latter case happens if $p_{n,\rho_{i}}[s]\neq p_{n,\rho_{i}}[s-1]$.)
In this case, $x_n$ \emph{abandons} this $y_{m'}$, if it existed at all,
and we will need to choose a new value $m$ for $h_{s+1}(n)$.
The element $y_{m'}$ becomes \emph{unattached},
and all lower-priority requirements will be injured at this stage.

If $p_{n,\rho_i}[s]$ is the zero polynomial, then $x_n$ currently appears
to be differentially transcendental, so we set $h_{s+1}(n)$
equal to the least number $m$ such that $y_m\notin F_s$.
Elements already in $F_s$ already satisfy a polynomial,
so we cannot define $h_{s+1}(n)$ to be an existing $m$.
The new $y_m$ is adjoined to $F_{s+1}$, with no change to $U_{s+1}$
(so that $y_m$ appears differentially transcendental in $F_s$).

If $p_{n,\rho_i}[s+1]$ was nonzero, then we wish to find some $y_m$ for
which we can make $p_{n,\rho_i}(y_{h(n_0)},\ldots,y_{h(n_{i-1})},Y)[s]$
the minimal differential polynomial over $\Q\{ y_{h(n_0)},\ldots,y_{h(n_{i-1})}\}[s]$.
For $p_{n,\rho_i}$ of positive order, this can be done by taking $m=r+1$ if needed,
since no facts about $y_{r+1}$ have yet been stated in $F$.  (Algebraic polynomials
$p_{n,\rho_i}$ will have no more roots in $K$ than they are allowed
to have in $F$, so either $y_{r+1}$ or an existing $y_m$ must suffice.)
However, for the sake of $\R_n$, we need to choose $m$
as small as possible without injuring higher-priority requirements.
It is now necessary to define the process by which $\R_n$ \emph{asks permission}
from those requirements to add a polynomial to $U_{s+1}$; this appears directly below.
For the least $m\leq r+1$ such that $\S_{n-1}$ (and hence all higher-priority requirements)
grant permission, and such that $y_m$ is not yet a root of any lower-order
polynomial than $p_{n,\rho_i}$, we adjoin
$$p_{n,\rho_{i}}(Y_{h(n_{0})},\ldots,Y_{h(n_{i-1})},Y_m)[s]$$
to $U_{s+1}$; this means we are setting
$p_{n,\rho_{i}}(y_{h(n_{0})},\ldots,y_{h(n_{i-1})},y_m)=0[s]$
in $F$, just as $p_{n,\rho_{i}}(x_{n_{0}},\ldots,x_{n_{i-1}},x_n)=0[s]$
in $K_s$.  With $h_{s+1}(n)=m$, our $h_{s+1}$ still
defines a partial isomorphism, based on the approximation $K_s$.
If $m=r+1$, we also add $x_{r+1}$ to $F_{s+1}$.

No matter which case held in item (3), we do not go on to the next substage,
but continue instead with the Final Step of stage $s+1$ (described below).

\end{enumerate}

This covers all the possibilities at substages dedicated to $\R$-requirements.
Notice that, even if $m$ lay in $\rg{h_s}$ but not in $\rg{h_{s+1}}$,
$y_m$ is still in $F_{s+1}$, and $U_s\subseteq U_{s+1}$.
This is necessary in order for $F$ to be computable.
Eventually, $\S_m$ will choose an $h$-preimage for $m$
respecting these conditions.

\textbf{Asking permission to adjoin to $U$.}
\label{permission}
Suppose $g\in\Q\{Y_{m_0},\ldots,Y_{m_i})$ is a polynomial which we wish
to add to $U_{s+1}$.  To \emph{ask permission} from a requirement $\R_n$
or $\S_m$ to do this, we choose the unique $i$ with $m_{i,s}=m$ (for $\S_m$)
or with $h_s(n)=m_{i,s}$ (for $\R_n$), and run the following process.
If $\DCF0\proves\psi_s\to g\neq 0$, then permission is immediately denied.  
Otherwise, let $E_{0,s}=\Q$, and define $E_{i,s}$ by recursion on $j < i$.
\begin{itemize}
\item
If $n_j=h^{-1}(m_j)\leq m_j[s]$,
then $\R_{n_j}$ controls $y_{m_j}$, and we set $E_{j+1,s}$ to be the computable
differential field $E_{j,s}\{ y_{m_j}\}/\la p_{n_j,\rho_j}\ra$,
whose atomic diagram $\Delta(E_{j+1,s})$ is generated over $\DCF0\cup\Delta(E_{j,s})$
by the formula $p_{n_j,\rho_j}(y_{m_0},\ldots,y_{m_j})=0$ along with
the statements that $y_{m_j}$ is not a zero of any polynomial over $E_{j,s}$
of lower order than this $p_{n_j,\rho_j}$.
\item
If $n_j=h^{-1}(m_j)> m_j[s]$,
then $\S_{m_j}$ controls $y_{m_j}$, and we set $E_{j+1,s}$ to be the
computable differential field extending $E_{j,s}$ with one new
generator $y_{m_j}$ satisfying the type given by $\gamma(e_{j,s},f)$,
where $\gamma$ is the type function from Lemma \ref{lemma:typefct},
$e_{j,s}$ is an index for $\Delta(E_{j,s})$, and $f$ is the current minimal differential
polynomial of $y_{m_j}$ in $F_s$, as given by Lemma \ref{lemma:minpoly}.
\end{itemize}
So $E_{i,s}$ is the differential field which the higher-priority requirements
currently believe us to be building.  (If the approximations given by $K$
subsequently change, though, then $E_{i,s}$ could turn out not to be a subfield
of the $F$ we finally build.)  Hence the theory $\DCF0\cup\Delta(E_{i,s})$ is complete
and consistent, is decidable uniformly in $i$ and $s$ using quantifier elimination
in \DCF0, and contains constant symbols
$y_{m_0},\ldots,y_{m_i}$.  Now $g(y_0,\ldots,y_r)$ may have more variables
than just these constants, so we check whether the formula
$$ (\exists y_{k_0}\exists y_{k_1}\cdots\exists y_{k_l}) [\psi_s~\&~g=0]$$
lies in this theory, where $\{ k_0,\ldots,k_l\}=\set{k\leq r}{k\notin\{m_0,\ldots,m_i\} }$.
If so, then the requirement \emph{allows} $g$ to be adjoined to $U_{s+1}$;
if not, then it \emph{denies permission} for this adjoinment.
This completes the process of asking permission.
(Notice that in fact we have received permission not just from
the given requirement $\R_n$ or $\S_m$, but from all
higher-priority requirements as well, via their subfields $E_{j,s}$ of $E_{i,s}$.)

\textbf{$\S_m$-substages.}
Next we explain the instructions for a substage for the requirement $\S_m$.
We fix the $i$ (which must exist) such that $m_{i,s}=m$, and the current
minimal differential polynomial $f$ of $y_m$ over  $y_{m_0},\ldots,y_{m_{i-1}}[s]$.
Now either $h_{s+1}^{-1}(m)$ has already been determined
by some higher-priority $\R_n$ (so $\S_m$ has nothing to do),
or $h_s(n)=m$ for some $n>m$, or $y_m$ is currently unattached
(i.e., $m\notin\rg{h_s}$).  In these latter two cases, it is not clear
that we will ever be able to find any $x\in K$ with minimal
differential polynomial $f$ over $x_{n_0},\ldots,x_{n_{i-1}}[s]$,
since $f$ might not be constrainable over these elements.
(If $h_s^{-1}(m)=n$ is defined, then $x_n$ currently appears
to fill this role, but in the noncomputable differential field $K$,
this could change at any time.)  So the requirement $\S_m$ 
will use the type function $\gamma$ from Lemma \ref{lemma:typefct},
knowing that $\gamma$ must give us
an index for a complete principal $1$-type over $E_{i,s}$
which is consistent with $\psi_s$ (and in particular with $f=0$).

At a substage for a requirement $\S_m$ within stage $s+1$,
we follow these instructions.  Fix the unique $i$ such that $m=m_{i,s}$.
If there exists an $n\leq m$ such that $h_{s+1}(n)$ has already been defined
to equal $m$, then we go on to the next substage.  Also,
if $h_s^{-1}(m)$ was defined and equal to some $n=n_{i,s}>m$,
and $p_{n,\rho_{i}}[s]\neq p_{n,\rho_{i}}[s-1]$, then $y_m$ becomes
\emph{unattached}.  We make $h_{s+1}^{-1}(m)$ undefined and end this substage,
and, instead of continuing to the next substage, we execute the Final Step
of stage $s+1$.

Otherwise we create the computable differential field $E_{i,s}$ currently envisioned
by the higher-priority requirements, exactly as defined above in the process
for asking permission from the next-higher-priority requirement
$\R_m$.  Let $e_{i,s}$ be an index for the atomic diagram $\Delta(E_{i,s})$.
For each of the first $s$ irreducible differential polynomials
$q_0,\ldots,q_s\in\Q\{ Y_{m_0},\ldots,Y_{m_i}\}$
of strictly lower order than $f$ in $Y_{m_i}$, we compute
$\phi_{\gamma(e_{i,s},f)}(\ulcorner \psi_s~\&~q_j=0\urcorner)$; that is, we ask whether
the formula $(\psi_s~\&~q_j=0)$ belongs to the $1$-type determined by $\gamma$
for $y_{m_i}$ over $E_{i,s}$, given that $f(y_{m_0},\ldots,y_{m_{i-1}},Y)$ is
currently the minimal differential polynomial of $y_{m_i}$ over $E_{i,s}$.
If so, then for the least such $j$, we adjoin $q_j$ to $U_{s+1}$,
having already seen from $E_{i,s}$ that this will not injure any higher-priority requirements;
we then end this substage and go directly to the Final Step
of the stage.  (This constitutes an injury to all lower-priority
requirements, but since the order of the minimal polynomial of $y_{m_i}$
can only decrease finitely often, there will be only finitely many such injuries.)

If there is no $j\leq s$ for which $(\psi_s~\&~q_j=0)$ belongs
to the $1$-type in question, then we keep $U_{s+1}=U_s$,
and act according to the following three cases, which together
complete the instructions for the $\S_m$-substage.
\begin{enumerate}
\item
If $h_s^{-1}(m)$ was defined and equal to some $n=n_{i,s}>m$,
and no $n'<n$ with $n'\notin\{ n_0,\ldots,n_{i-1}\}[s]$ has $p_{n',\rho_i}[s]$ equal
to the apparent minimal differential polynomial $f$ of $y_m$ over
$\{y_0,\ldots,y_{m_{i-1}}\}$ in $F_s$, then we keep $h_{s+1}(n)=m$
and go on to the next substage.
\item
If $h_s^{-1}(m)$ was defined and equal to some $n=n_{i,s}>m$,
and some $n'<n$ with $n'\notin\{ n_0,\ldots,n_{i-1}\}[s]$ has $p_{n',\rho_i}$ equal
to the apparent minimal differential polynomial $f$ of $y_m$ over
$\{y_0,\ldots,y_{m_{i-1}}\}$ in $F_s$, then $y_m$ becomes \emph{unattached}.
We make $h_{s+1}^{-1}(m)$ undefined and end this substage,
and, instead of continuing to the next substage, we execute the Final Step
of stage $s+1$.
(At the $\S_m$-substage of the next stage, we will search
for a new $h$-preimage for $m$, most likely the $n'$ found above.)

\item
Otherwise, $h_s^{-1}(m)$ was undefined and $h_{s+1}^{-1}(m)$
has not been defined at an earlier substage of this stage.
We check to see whether any $n\leq s$ with $n\notin\{ n_0,\ldots,n_{i-1}\}[s]$
has $p_{n,\rho_{i}}[s]=f$.  If so, then we define $h_{s+1}(n)=m$
(for the least such $n$); if not, then $h_{s+1}^{-1}(m)$ remains
undefined.  In either case we proceed to the Final Step.
(Eventually some such $n$ will have to reveal itself, since,
once our choice of $f$ has stabilized, this $f$ will be
constrainable over the higher-priority elements of $F$,
hence must have a zero in $K$ over the corresponding elements there.)
\end{enumerate}

\textbf{Final Step.}
To finish stage $s+1$, after completing the last substage,
consider the next differential polynomial $g_s(Y_0,\ldots,Y_k)$
in a fixed computable enumeration $g_0,g_1,\ldots$ of
$\Q\{ Y_0,Y_1,\ldots\}$.  Consider the lowest-priority
element $y_{r'}$ currently in $F_s$.
We ask permission either from the requirement $\R_n$
(where $h_{s+1}(n)=r'\geq n$, if such an $n$ exists), or else
from the requirement $\S_{r'}$, to adjoin $g_s$ to $U$.
If this permission is granted, then $g_s\in U_{s+1}$.
If not, then $U_{s+1}$ stays unchanged and we know $g_s\notin U$.
(Thus $U$ will be decidable.)
This completes the Final Step, and ends stage $s+1$.

We set $F=\set{y_m}{m\in\omega}$,
but the important objects constructed were the decidable set $U=\bigcup_s U_s$
and the finite functions $h_s$, whose limit will be the isomorphism from $K$ onto $F$.
Notice that, every time any differential polynomial $g(Y_0,\ldots,Y_{k_f})$
was enumerated into $U_s$, the permission process confirmed that the formula
$$ \exists Y_0\cdots\exists Y_r  (\psi_s\wedge g=0)
$$
belonged to the theory \DCF0.
It follows that the entire set of formulas $\psi_s$, for all $s$, is consistent with \DCF0.

The bijection between $F$ and $K$
will follow once we prove these claims for all $i$:
\begin{itemize}
\item
$n_i=\lim_s n_{i,s}$ exists, and the map $i\mapsto n_i$
is a permutation of $\omega$;
\item
$m_i=\lim_s m_{i,s}$ exists, and the map $i\mapsto m_i$
is a permutation of $\omega$;
\item
the function $h=\lim_sh_s$ is a bijection from 
$\omega$ onto $\omega$, and hence defines a bijection $x_n\mapsto y_{h(n)}$
from $K$ onto $F$; and
\item
the limit $p_i=\lim_s p_{n_{i,s},\rho_{i,s},s}\in\Q\{ X_{n_0},X_{n_1},\ldots,X_{n_i}\}$ exists,
and 
$U$ contains $p_i(Y_{h(n_0)},\ldots,Y_{h(n_i)})$, and no $q(Y_{h(n_0)},\ldots,Y_{h(n_i)})$
in $U$ has lower $Y_{h(n_i)}$-rank than $p_i$.
(Here $\rho_i=(n_0,\ldots,n_{i-1})=\lim_s\rho_{i,s}$, from the first claim.)
\end{itemize}
The first three claims here can be proven together by a single induction.

\begin{lemma}
\label{lemma:reqts}
For every $m$, there exists a unique $i$ with $\lim_s m_{i,s}=m$;
likewise, for every $n$, there exists a unique $i$ with $\lim_s n_{i,s}=n$.
Thus every requirement $\R_n$ and $\S_m$ is satisfied
by the foregoing construction.
\end{lemma}
\begin{pf}
The uniqueness of $i$, for any single $m$ or $n$, is immediate
from our definitions of $m_{i,s}$ and $n_{i,s}$.  We specifically
excluded all repetitions from the first sequence, making $m_{i,s}\neq m_{j,s}$
for every $i<j$, and we made every $h_s$ injective.
Recall that by our definition, at stage $s$, every $n_{i,s}$
except the very last one lies in $\dom{h_s}$.
The injectivity of each $h_s$ follows from its construction:  we always
included in $\psi_s(Y_0,\ldots,Y_r)$ the conditions that $Y_i\neq Y_j$
for all $i<j\leq r$, and similarly in $\sigma_{i,s+1}$
that $X_{n_{i,s}}\neq X_{n_{j,s}}$,
and then we required the choice of each new $h_{s+1}(n)$ to have
$\sigma_{i,s+1}(Y_{h_{s+1}(n_{0,s})},\ldots,Y_{h_{s+1}(n_{i-1,s})},Y_{h_{s+1}(n)})$
consistent with $\psi_s(Y_0,\ldots,Y_r)$.

We proceed by induction on these requirements, according to their priority order,
starting with $\R_0$.  The inductive hypothesis is that there exists a stage $s_0$
such that, for every $s\geq s_0$ and each higher-priority requirement $\R_{n'}$ or $\S_{m'}$,
there are unique numbers $j$ and $k$ with $n_{j,s}=n'$ and $m_{k,s}=m'$ and
$h_s(n')=h_{s_0}(n')$ and $h_s^{-1}(m')=h_{s_0}^{-1}(m')$.
Turning to the minimal polynomials in $K$,
we may also assume that $s_0$ is so large that, for every $n'=n_{j,s}<n$,
$p_{n',\rho_{j,s},s}=p_{n',\rho_{j,s},s_0}$ (noting that $\rho_{j,s}=\rho_{j,s_0}$
by the previous part of the hypothesis).  That is, all approximations
to minimal polynomials of higher-priority elements of $K$ have converged
by stage $s_0$.  It follows that, from stage $s_0+1$ on, every substage
for a higher-priority requirement will do nothing.  Moreover, at all subsequent stages $s$,
the field $E_{i,s}$ will have stabilized as one particular differential subfield $E_i$ of $F$
(where $i$ is chosen so that either $m=m_{i,s}$ or $h_s(n)=m_{i,s}$).

Suppose this inductive hypothesis holds of every requirement of higher
priority than $\R_n$.  If there exists an $m<n$ with $h_{s_0}(n)=m$,
then the satisfaction of $\S_m$ shows that $\R_n$ is satisfied as well.  So assume
that there is no such $m$.  
Let $\rho=\rho_{i,s_0+1}$ be the sequence of indices of elements in $K$
of higher priority than $n$.  This too never changes at stages $>s_0$.
But now the approximations $p_{n,\rho,s}$ to the minimal differential polynomial
of $x_n$ over $\Q\la x_0,\ldots,x_{i-1}\ra$ (with $x_j=\lim_s x_{j,s}$) must converge,
to some limit $p_n(X_0,\ldots,X_i)$.  Let $s_1>s_0$ be a stage by which this
convergence has occurred.  
If $h_{s_1}(n)$ is undefined, then at stage $s_1+1$
the construction will reach the substage for $\R_n$ and will act according
to item (3) at that substage, and will choose a value $h_{s_1+1}(n)\leq r+1$.
This $y_{h_{s_1+1}(n)}$ therefore lies in $F_s$ at all $s\geq s_1+1$.
At the next stage $s_1+2$, $n$ will lie in the domain of $h_{s_1+1}$,
and therefore will have $n=n_{i,s_1+1}$ for some $i$, i.e., $n$ will have been
assigned a priority, corresponding to the requirement $\R_n$.
From then on, item (2) in the substage for $\R_n$
will always apply, leaving the value of $h_s(n)$ unchanged.
Moreover, in the process of asking permission, $E_{i,s}$ ensures that
the minimal polynomial of $y_{h_s(n)}$ in $F$ would only change
if the rank of a higher-priority element changed, or if the approximation
to $p_n$ changed.  By assumption neither of these ever changes again,
so the minimal polynomial of $y_{h_s(n)}$ in $F$ stays fixed forever.
Therefore, $h_s(n)$ will never again change its value, and
the requirement $\R_n$ is indeed satisfied.  The existence
of the (unique) $i$ with $n=n_i=\lim_s n_{i,s}$ follows.

Now we turn to the inductive step for a requirement $\S_m$,
using the stage $s_0$ defined above by the inductive hypothesis
on all higher-priority requirements.  Once again, it follows that every
higher-priority requirement will do nothing at its substage during
each stage $>s_0$, and so the $\S_m$-substage will be reached
at every such stage.  If $h_{s_0}(n)=m$ for some $n\leq m$,
then the satisfaction of the higher-priority requirement $\R_n$
shows that $m=\lim_s h_s(n)$; so assume that this is not the case.
Now $F_s$ increases at infinitely many stages $s$,
so eventually some $F_{s_1}$
will include $y_m$.  At this point, an $i$ will be chosen
for which $m_{i,s}=m$, since this happens for all indices
of elements of $F_s$.  Moreover, taking $s_1>s_0$
and knowing that the higher-priority
requirements never act again, we will have $m_{i,s}=m$
at all stages $>s_1$ as well; this proves existence of
the $i$ with $m=m_i=\lim_s m_{i,s}$, and its uniqueness
was already seen.

At stage $s_1$, $y_m$ has an apparent minimal differential
polynomial $f$ over the higher-priority requirements.  Since
$E_i$ never again changes, and every subsequent adjoinment to $U$
will require the permission of $\S_m$, we know that $y_m$
must realize the type $\Gamma$ over $E_i$ given by the type function:
$\phi_{\gamma(e_i,f)}$ computes this type.  Since $\Gamma$ is principal,
there must exist an $s$ such that $\Gamma$ contains a formula 
of the form $(q=0~\&~\psi_s)$ which generates $\Gamma$.
This $q$ is therefore constrainable (with $\psi_s$ providing
the constraint, if a nontrivial one is needed), and when $q$ appears
in an $\S_m$-substage, $y_m$ will be defined to be a zero of this $q$.

(Lemma \ref{lemma:typefct} did not actually claim that, whenever
$q=0$ lies in the type $\gamma(e,f)$ with $q$ of smaller rank than $f$,
the index $\gamma(e,q)$ must then define the same type as
$\gamma(e,f)$.  It can readily be arranged for this to be so, however;
and even if it were not so, it would only contribute finitely many more
injuries to the lower-priority requirements.)

So eventually $y_m$ is found to be a zero of a constrainable $q$,
in particular, of the smallest-rank $q$ such that $q=0$ lies in this type.
Once this has happened, the differentially closed field $K$ must reveal an $x_n$
realizing this same type over the $h$-preimages of the higher-priority elements.
For the least such $n$, once the $K$-approximation settles on $q$
as the minimal differential polynomial of this $x_n$ (and once all $x_{n'}$
with $n'<n$ have settled on their own minimal differential polynomials
distinct from $q$), we will define $h(n)=m$,
and will preserve $h(n)=m$ forever after.  This completes the proof of the lemma.
\qed\end{pf}

Finally we consider the last claim, for a fixed $i$.
The first part of the claim has already been noted:
we have seen above that the limit $n_j=\lim_s n_{j,s}$
exists for every $j$, and so, with $\rho_i=\lim_s\rho_{i,s}$,
the computable approximations in $K$ all converge to the actual
minimal differential polynomials $p_i=\lim_s p_{n_i,\rho_{i,s},s}$.
We have also seen above that $m_i=h(n_i)=\lim_s h_s(n_i)$ exists.
But each map $h_s$ defines a partial isomorphism
from the approximation $K_s$ into $F$, and so,
once all the approximations for a given fragment of $K$
have converged, the limit $h$ on this fragment will define
a partial isomorphism.  Since $h$ is also a bijection,
it does in fact define an isomorphism $x_n\mapsto y_{h(n)}$.
\comment{
If $m_i<n_i$, then the reason why $\S_{m_i}$ chose to make
and keep $h_s^{-1}(m_i)=n_i$ is that $p_i(y_{m_0},\ldots,y_{m_i})\in U$
(as seen in the instructions for unattached $y_m$, where $n_i$
was first chosen) and that item (3) of the $\S_{m_i}$-substage
never found a differential polynomial of lower rank for which
$y_{m_i}$ could be a zero.  So there is no such polynomial,
which is exactly what the claim requires.

On the other hand, if $m_i\geq n_i$, then it was $\R_{n_i}$
which chose to make and keep $h_s(n_i)=m_i$.
It did this at $\R_{n_i}$-substages because $p_i$ appeared to be the minimal
polynomial of $y_{m_i}$ over $y_{m_0},\ldots,y_{m_{i-1}}$.
Therefore, the construction never made $y_{m_i}$ a zero
of any lower-rank differential polynomial.
}
This completes the proof of the final claim.

It follows that the operations in $F$ are computable.  For instance,
given any $y_i,y_j\in F$, the elements $x_{h^{-1}(i)}$ and $x_{h^{-1}(j)}$
of $K$ have a sum $x_k$.  Since $h$ defines an isomorphism,
the polynomial $Y_i+Y_j-Y_{h(k)}$ must lie in the decidable set $U$,
and when we find it, we will know that $y_{h(k)}=y_i+y_j$.
Multiplication and differentiation are similarly computable, so
$F$ is a computable structure, and the isomorphism $h$ from $K$
onto $F$ establishes Theorem \ref{thm:low}.
\qed\end{pf}

Theorem \ref{thm:low} will remind many readers of the well-known
theorem of Downey and Jockusch from \cite{DJ94},
that every low Boolean algebra has a computable copy.
However, the parallels between these results are few.  The latter theorem
has been extended to included low$_4$ Boolean algebras,
in work by Thurber \cite{T95} and Knight and Stob \cite{KS00}, whereas by
Theorem \ref{thm:jumpspectra}, the result for \DCF0 does not
even extend to the low$_2$ case.  Moreover, the proof of Theorem \ref{thm:low}
constructed a $\Delta^0_2$-isomorphism from the low model of \DCF0 to
its computable copy, whereas for Boolean algebras, there is always a
$\Delta^0_3$-isomorphism but not always a $\Delta^0_2$ one.
The construction here relied heavily on the completeness
and decidability of the theory \DCF0, whereas the theory of Boolean algebras
is certainly not complete.  Conversely, the construction in \cite{DJ94} uses theorems
of Vaught and Remmel which are specific to Boolean algebras,
with no obvious analogue for \DCF0.

The closer analogy is to the theory \ACF0,
for which Theorem \ref{thm:low} is trivially true,
since \emph{every} countable algebraically closed field
has a computable presentation.  All those of finite transcendence
degree over $\Q$ are relatively computably categorical,
meaning that every presentation of degree $\bfd$ has 
 a $\bfd$-computable isomorphism onto a computable copy.  
The unique countable model of \ACF0 of infinite transcendence degree
over $\Q$ is not, but it is relatively
$\Delta^0_2$-categorical, since in one jump over the atomic diagram of the structure,
one can compute a transcendence basis for the field over $\Q$.
For low models of \ACF0, one can give a much simpler version of the
priority construction used in Theorem \ref{thm:low}.
For readers who find the construction in the proof of Theorem \ref{thm:low}
daunting, carrying out this construction for \ACF0 might be a useful prelude.

\section{Spectra of Differentially Closed Fields}
\label{sec:conclusion}

\begin{prop}
\label{prop:relative}
For every countable model $K$ of \DCF0 of Turing degree $\bfc$,
every degree $\bfd$ with $\bfd'\geq\bfc'$ lies in the spectrum of $K$.
\end{prop}
\begin{pf}
One simply runs the same construction as in Theorem \ref{thm:low},
relative to an oracle from $\bfd$.  Since $\bfd'\geq\bfc'$, this oracle can
compute all the necessary approximations to facts about $K$ and about
minimal differential polynomials in $K$, so this produces a
$\bfd$-computable differential field isomorphic to $K$.  As mentioned
in Subsection \ref{subsec:CMT}, Knight's theorem from \cite{K86}
then shows that $\bfd\in\spec{K}$, since no differentially closed
field is automorphically trivial.
\qed\end{pf}
\begin{defn}
\label{defn:FJE}
\emph{First-jump equivalence} is the relation 
$\sim_1$ on Turing degrees:
$$ \bfc\sim_1\bfd\iff\bfc'=\bfd'.$$
\end{defn}
Proposition \ref{prop:relative} shows that every spectrum of a model
$K$ of \DCF0 \emph{respects} $\sim_1$, in the sense that,
whenever $\bfc\sim_1\bfd$, we have $(\bfc\in\S\iff\bfd\in\S)$.  It follows that
$\spec{K}$ is actually determined by its \emph{jump spectrum}
$\set{\bfd'}{\bfd\in\spec{K}}$.  Moreover, this proposition, along with
Lemma \ref{lemma:oracle} (which is easily proven using the methods
of \cite[Chapter VI]{S87}), yields a quick proof
of a property for \DCF0 which was already known to hold for linear
orders, Boolean algebras, and trees (viewed as partial orders), by results of
Richter in \cite{R81}.  When the question of spectra of differentially closed fields
first arose, this corollary was quickly observed by Andrews
and Montalb\'an, who pointed out that it follows from \cite{R81}.

\begin{cor}[cf.\ Andrews \& Montalb\'an]
\label{cor:nodegree}
No countable differentially closed field $K$ of characteristic $0$
intrinsically computes any noncomputable set $B\subseteq\omega$.
That is, the spectrum of $K$ cannot be contained within the upper
cone $\set{\bfd}{\bfb\leq\bfd}$ above a nonzero degree $\bfb$.  
In particular, if such a spectrum
has a least degree under $\leq_T$ among its elements, then that degree is $\bfz$.
\end{cor}
\begin{pf}
Let $K$ have degree $\bfc$.
Lemma \ref{lemma:oracle} below yields a degree
$\bfd$ with $\bfb\not\leq\bfd$ and $\bfc'\leq\bfd'$.
But then $\bfd\in\spec{K}$ by Proposition \ref{prop:relative}.
\qed\end{pf}
\begin{lemma}[Folklore]
\label{lemma:oracle}
For every noncomputable set $B$ and every set $C$, there exists some set $D$
with $B\not\leq_T D$ and $C'\leq_T D'$.  Indeed $C'\leq_T\emptyset'\oplus D$.
\qed\end{lemma}
The main consequence of Proposition \ref{prop:relative} is a very
precise description of the spectra of models of \DCF0 in terms of arbitrary spectra.
Theorem \ref{thm:HKSS} shows that items (2) and (3) of Theorem \ref{thm:allspectra}
could equally well allow $G$ and $J$ to vary over structures in all computable languages.
\begin{thm}
\label{thm:allspectra}
For a set $\S$ of Turing degrees, the following are equivalent.
\begin{enumerate}
\item
$\S$ is the spectrum of some countable model $K$ of \DCF0.
\item
There exists a countable, automorphically nontrivial graph $G$ for which
$ \S=\set{\bfd}{\bfd'\in\spec{G}}$.
\item
$\S$ respects $\sim_1$ and there exists a countable,
automorphically nontrivial graph $J$ with $\S=\spec{J}$.
\end{enumerate}
\comment{
The spectra of countable differentially closed fields of characteristic $0$
are precisely the sets of the form:
$$ \set{\bfd}{\bfd'\in\spec{G}},$$
where $G$ ranges over all countable, automorphically nontrivial
undirected graphs (or equivalently, over all countable,
automorphically nontrivial structures in computable languages).
}
\end{thm}
\begin{pf}
The implication $(2)\!\!\implies\!\!(1)$ is precisely Theorem \ref{thm:jumpspectra}
above.  Also, $(1)\!\!\implies\!\!(3)$ follows from Proposition \ref{prop:relative}
and Theorem \ref{thm:HKSS}.  To establish $(3)\!\!\implies\!\!(2)$,
given $J$, we appeal to the following theorem, proven by Soskova and Soskov in \cite{SS09}
and independently by Montalb\'an in \cite{M09} and first presented by Soskov
in a talk in 2002. 
\begin{thm}[see \cite{M09,SS09}]
\label{thm:jumps}
For every countable structure $\A$, there exists a countable
structure $\A'$, the \emph{jump of the structure $\A$}, such that
$\spec{\A'}=\set{\bfc'}{\bfc\in\spec{\A}}$.
\end{thm}
Using Theorem \ref{thm:HKSS}, we convert the jump $J'$
of our $J$ into a graph $G$, with $\spec{G}=\set{\bfc'}{\bfc\in\spec{J}}$.
Since $J$ is automorphically nontrivial, so is $G$.
Now each $\bfd\in\S=\spec{J}$ has $\bfd'\in\spec{G}$.  Conversely,
for every $\bfd$ with $\bfd'\in\spec{G}$, we have some
$\bfc\in\spec{J}=\S$ with $\bfc'=\bfd'$, making $\bfd\in\S$
since $\S$ respects $\sim_1$.
\qed\end{pf}

\end{document}